\documentclass[12pt,reqno]{amsart}
\usepackage[utf8]{inputenc}
\usepackage{fancyhdr}
\usepackage{amssymb}
\usepackage{amsthm}
\usepackage{amsmath}
\usepackage{mathtools}
\usepackage{amsxtra}
\usepackage{mathrsfs}
\usepackage[all]{xy}
\usepackage{colonequals}
\usepackage{enumitem}
\usepackage{fullpage}
\usepackage{breqn}
\usepackage{caption}
\usepackage{placeins}

\theoremstyle{plain}
\newtheorem{thm}{Theorem}[section]
\newtheorem{lem}[thm]{Lemma} 
\newtheorem{prop}[thm]{Proposition}
\newtheorem{cor}[thm]{Corollary}
\newtheorem{conj}{Conjecture}[section]
\theoremstyle{definition}
\newtheorem*{ack*}{Acknowledgment}
\theoremstyle{remark}

\DeclareMathOperator{\id}{id}

\DeclareMathOperator{\li}{li}

\newcommand{\N}{\mathbb N}
\newcommand{\PP}{\mathbb P}

\newcommand{\R}{\mathbb R}
\newcommand{\Z}{\mathbb Z}

\newcommand{\rad}{\text{rad}}

\newcommand{\psmod}[1]{\,(\textup{\text{mod}}\,{#1})}

\input cyracc.def

\usepackage[hidelinks]{hyperref}
\makeatletter
\@namedef{subjclassname@2020}{\textup{2020} Mathematics Subject Classification}
\makeatother

\begin{document}
\title{The shifted prime-divisor function over shifted primes}
\author{Kai (Steve) Fan}
\address{Max-Planck-Institut f\"{u}r Mathematik, 53111 Bonn, Germany}
\email{steve.fan@mpim-bonn.mpg.de}
\setlength{\footskip}{6mm}
\begin{abstract}
Let $a,b\in\Z\setminus\{0\}$. For every $n\in\N$, denote  by $\omega_a^*(n)$ the number of shifted-prime divisors $p-a$ of $n$, where $p>a$ is prime. In this paper, we study the moments of $\omega_a^*$ over shifted primes $p-b$. Specifically, we prove an asymptotic formula for the first moment and upper and lower bounds of the correct order of magnitude for the second moment. These results suggest that the average behavior of $\omega^*_a$ on shifted primes is similar to its average behavior on natural numbers. We shall also prove upper bounds for the mean values of sub-multiplicative functions in a nice class over the least common multiples of the shifted primes $p-a$ and $q-b$. Such upper bounds are intimately related to the second moments of $\omega^*_a$ over natural numbers and over shifted primes. Finally, we propose a new conjecture on the second moment of $\omega_1^*$ over natural numbers and provide a heuristic argument in support of this conjecture.
\end{abstract}

\maketitle
\section{Introduction}\label{S:Intro}
For every $n\in\N$, let $\omega^*(n)$ denote the number of
shifted prime-divisors $p-1$ of $n$, that is, 
\[\omega^*(n)\colonequals\sum_{(p-1)\mid n}1.\]
Equivalently, $\omega^*(n)$ counts the number of positive divisors $d\mid n$ such that $d+1$ is prime. This function has interesting applications in primality testing \cite{APR}, and a variation of it which counts the number of positive divisors $d\mid n$ such that $dk+1$ is prime for a given $k\in\N$ has been proven useful in the study of the distribution of Carmichael numbers \cite{AGP}. Although its definition bears a resemblance to that of the prime divisor function $\omega(n)$, which counts the number of distinct prime factors of $n$, the behavior of the function $\omega^*(n)$ is in some sense closer to that of the divisor function $\tau(n)$, which counts the number of all the positive divisors of $n$. This is especially the case when it comes to maximal orders. Indeed, Prachar \cite{Pra} showed that for infinitely many $n$, we have
\begin{align*}
	\omega^*(n)&>\exp\left(c_1\frac{\log n}{(\log\log n)^2}\right) \text{~~~(unconditionally)},\\
	\omega^*(n)&>\exp\left((\log\sqrt{2}-\epsilon)\frac{\log n}{\log\log n}\right) \text{~~~(under GRH)},
\end{align*}
where $c_1>0$ is some absolute constant, and $\epsilon>0$ is fixed but otherwise arbitrary. Prachar's unconditional lower bound was later improved by Adleman, Pomerance and Rumely \cite[Proposition 10]{APR} to 
\[\omega^*(n)>\exp\left(c_2\frac{\log n}{\log\log n}\right) \]
for infinitely many $n$, where $c_2>0$ is some absolute constant, matching Prachar's conditional lower bound under GRH except for the unknown value of $c_2$. On the other hand, it is well known that 
\begin{align*}
	\limsup\limits_{x\to\infty}\frac{\omega(n)}{\log n/\log\log n}&=1,\\
	\limsup\limits_{x\to\infty}\frac{\log\tau(n)}{\log n/\log\log n}&=\log 2,
\end{align*}
the second of which is due to Wigert in 1907. Thus, the maximal order of $\log\omega^*(n)$ is comparable to those of $\omega(n)$ and $\log\tau(n)$. 

Another aspect in which $\omega^*$ and $\tau$ share similarities concerns moments. For every $k\in\N$, we define the {\it $k$th moment} of $\omega^*$ over natural numbers by
\[M_k(x)\colonequals\frac{1}{x}\sum_{n\le x}\omega^*(n)^k.\]
Like the first moment of $\omega$, the first moment of $\omega^*$ is asymptotically $\log\log x$ \cite{Pra}. Indeed, a quick application of Mertens' second theorem shows $M_1(x)=\log\log x+O(1)$. However, the second moment of $\omega^*$ turns out to be much larger than that of $\omega$. The study of $M_2(x)$ was carried out in \cite{Pra,MM,Din} and culminated in the estimate $M_2(x)\asymp\log x$. One may compare this estimate with 
\[\frac{1}{x}\sum_{n\le x}\tau(n)^2\sim\frac{1}{\pi^2}(\log x)^3,\]
an old result due to Ramanujan. If we start with the asymptotic above and take into account the fact that there are two primality constraints in the expression of $\omega^*(n)^2$ given by 
\[\omega^*(n)^2=\sum_{\substack{p-1\mid n\\ q-1\mid n}}1,\]
then it is natural to predict that the order of magnitude of $M_2(x)$ is $\log x$, since the density of primes in $[1,x]$ is about $1/\log x$. Murty and Murty \cite{MM} conjectured that there exists a constant $C>0$ such that $M_2(x)\sim C\log x$. In \cite{DGZ} Ding, Guo, and Zhang argued that this conjecture holds with $C=2\zeta(2)\zeta(3)/\zeta(6)=315\zeta(3)/\pi^4$ based on the Elliott--Halberstam conjecture, where $\zeta$ is the Riemann zeta function. Unfortunately, it has been pointed out \cite{FP} that there is an error in \cite[Equation (4.8)]{MM} which is essentially the starting point of the argument in \cite{DGZ}. Due to this error, the constant obtained by Ding, Guo, and Zhang is probably incorrect, a conclusion which is also corroborated by the numerical data presented in \cite{FP}. In the same paper \cite{FP}, Pomerance and the author studied the third moment of $\omega^*$, proving that $M_3(x)\asymp(\log x)^4$. They also conjectured that for each integer $k\ge2$ one has $M_k(x)\sim C_k(\log x)^{2^k-k-1}$ for some constant $C_k>0$.
\par The material in the rest of the paper is organized as follows. In Section \ref{S:omega_a*(p-b)}, we study the moments of $\omega^*_a$ over shifted primes $p-b$, where $a,b\in\Z\setminus\{0\}$ are fixed, and \[\omega_a^*(n)\colonequals\sum_{\substack{(p-a)\mid n\\p>a}}1\]
counts the number of shifted-prime divisors $p-a$ of $n$, with $\omega_1^*=\omega^*$. Specifically, we shall prove an asymptotic formula for the first moment and upper and lower bounds of the correct order of magnitude for the second moment. These results suggest that the average behavior of $\omega^*_a$ on shifted primes is similar to its average behavior on natural numbers. 
\par Intimately related to the second moments of $\omega^*_a$ over natural numbers and over shifted primes are mean values of arithmetic functions over the least common multiples of shifted primes $p-a$ and $q-b$, where $a,b\in\Z\setminus\{0\}$. In Section \ref{S:f([p-a,q-b])}, we prove upper bounds for the mean values of sub-multiplicative functions from a nice class.
\par In Section \ref{S:2Momega*}, we return to the second moment $M_2(x)$ of $\omega^*$ over natural numbers which has been extensively studied in the literature \cite{Pra,MM,Din,DGZ,FP}. As mentioned above, the constant $C=2\zeta(2)\zeta(3)/\zeta(6)\approx3.88719$ obtained by Ding, Guo, and Zhang \cite{DGZ} is probably incorrect. Here we propose a new heuristic argument leading to the conjecture $C=\zeta(2)^2\zeta(3)/\zeta(6)=105\zeta(3)/4\pi^2\approx3.19709$, which is supported by the numerical computations conducted in \cite{FP}.
\par Finally, we conclude this paper with a brief discussion on the level sets of $\omega^*$ over shifted primes $p-b$ and possible analogues of the results obtained in \cite[Section 3]{FP}.\\
\\
\noindent {\bf Notation.} Throughout the paper, the symbols $\PP$, $\N$, $\Z$ and $\R$ stand for the set of prime numbers, the set of positive integers, the set of integers and the set of real numbers, respectively. For any $x\in\R$, we denote the integer part of $x$ by $\lfloor x\rfloor$, which is the largest integer not exceeding $x$, and the fractional part of $x$ by $\{x\}\colonequals x-\lfloor x\rfloor$. The letters $p,q,r$ always represent prime numbers, and $\pi(x)$ denotes the number of primes $p\le x$. For any $k\in\N$ and $a\in\Z$, the function $\pi(x;k,a)$ counts the number of primes $p\le x$ with $p\equiv a\psmod{k}$. For $m,n\in\Z\setminus\{0\}$, we write $\gcd(m,n)$ or simply $(m,n)$ when no confusion arises, for the greatest common divisor of $m$ and $n$, and $[m,n]$ for their least common multiple. For any $n\in\N$, we denote by $\rad(n)$ the radical of $n$, i.e., $\rad(n)\colonequals\prod_{p\mid n}p$. Given $p$, $\nu\in\N$ and $n\in\Z\setminus\{0\}$, the relation $p^{\nu}\parallel n$ means that $p^{\nu}\mid n$ whereas $p^{\nu+1}\nmid n$. We shall also make use of the identity function $\id(n)\colonequals n$, the M\"{o}bius function $\mu$, Euler's totient function $\varphi$, the sum-of-divisors function $\sigma$, and the $\kappa$-fold divisor function $\tau_{\kappa}$ for $\kappa\ge0$ with the abbreviation $\tau\colonequals\tau_2$. Besides, we denote by $\omega(n)$ the number of distinct prime factors of $n$ and by $\Omega(n)$ the total number of prime factors of $n$ counted with multiplicity. In addition, we write $1_{A}$ for the characteristic function of the set or condition $A$, meaning that $1_{A}(n)=1$ if $n$ belongs to the set $A$ or $n$ satisfies condition $A$ and $1_{A}(n)=0$ otherwise. Finally, we shall not only use Landau's big-$O$ notation and Vinogradov's notation $\ll$ interchangeably but also adopt the standard order notations $o,\gg,\asymp,\sim$ from analytic number theory. 
\medskip
\section{Moments of the shifted-prime divisor function over shifted primes}\label{S:omega_a*(p-b)}
In this section, we investigate the first and second moments of $\omega_a^*$ over shifted primes $p-b$, where $a,b\in\Z\setminus\{0\}$ are fixed. We begin with the following result which provides an asymptotic formula for the mean value of $\id/\varphi$ over shifted primes $p-a$ with the constraint $\gcd(p-a,b)=1$. The case $b=1$ is probably well-known. 
\begin{prop}\label{prop:shiftphi}
For any fixed $A_0>0$ and $a,b\in\Z\setminus\{0\}$, we have
\[\sum_{\substack{a<p\le x\\(p-a,b)=1}}\frac{p-a}{\varphi(p-a)}=C_{a,b}\li(x)+O\left(\frac{x}{(\log x)^{A_0}}\right)\]
for all $x\ge\max(a,2)$, where
\[\li(y)\colonequals\int_{0}^{y}\frac{dt}{\log t}\]
is the logarithmic integral, and 
\begin{equation}\label{eq:C_ab}
C_{a,b}\colonequals\prod_{p\nmid a}\left(1+\frac{1}{(p-1)^2}\right)\prod_{\substack{p\nmid a\\p\mid b}}\frac{(p-1)(p-2)}{(p-1)^2+1}.
\end{equation}
\end{prop}
\begin{proof}
Without loss of generality, we may assume that $b$ is square-free and that $x$ is sufficiently large. Let $x_a\colonequals x-a$. We start with the identity
\[\frac{n}{\varphi(n)}=\sum_{d\mid n}\frac{\mu(d)^2}{\varphi(d)}.\]
From this identity it follows that 
\begin{equation}\label{eq:sum id/phi}
\sum_{\substack{a<p\le x\\(p-a,b)=1}}\frac{p-a}{\varphi(p-a)}=\sum_{\substack{a<p\le x\\(p-a,b)=1}}\sum_{d\mid p-a}\frac{\mu(d)^2}{\varphi(d)}=S_1(x;a)+S_2(x;a)+O(\log x),
\end{equation}
where
\begin{align*}
S_1(x;a,b)&\colonequals\sum_{\substack{a<p\le x\\(p-a,b)=1}}\sum_{\substack{d\mid p-a\\d<\sqrt{p-a}}}\frac{\mu(d)^2}{\varphi(d)},\\
S_2(x;a,b)&\colonequals\sum_{\substack{a<p\le x\\(p-a,b)=1}}\sum_{\substack{d\mid p-a\\d>\sqrt{p-a}}}\frac{\mu(d)^2}{\varphi(d)}.
\end{align*}
Here we have used the estimate
\[\sum_{\substack{a<p\le x\\(p-a,b)=1\\ p-a=\square}}\frac{\mu(\sqrt{p-a})^2}{\varphi(\sqrt{p-a})}\le\sum_{n\le \sqrt{x_a}}\frac{\mu(n)^2}{\varphi(n)}\ll\log x,\]
The contribution from $S_2(x;a)$ is easily seen to be negligible. Indeed, we have
\begin{align*}
S_2(x;a,b)&=\sum_{\substack{a<p\le x\\(p-a,b)=1}}\sum_{\substack{d\mid p-a\\d<\sqrt{p-a}}}\frac{\mu((p-a)/d)^2}{\varphi((p-a)/d)}\\
&\ll\sum_{d\le\sqrt{x_a}}d\sum_{\substack{d^2+a< n\le x\\n\equiv a\psmod{d}}}\frac{\log\log(3(n-a))}{n-a}\\
&\ll \sqrt{x}(\log x)\log\log x,
\end{align*}
since $\varphi(n)\gg n/\log\log3n$.
\par It remains to estimate $S_1(x;a,b)$. To this end, we show that 
\begin{equation}\label{eq:S_101}
S_1(x;a,b)=C_{a,b}\li(x)+O\left(\frac{x}{(\log x)^{A_0}}\right),
\end{equation}
where $A_0>0$ is fixed but otherwise arbitrary. Note that
\[S_1(x;a,b)=\sum_{\substack{d<\sqrt{x_a}\\(d,a)=1}}\frac{\mu(d)^2}{\varphi(d)}\sum_{\substack{d^2+a<p\le x\\p\equiv a\psmod{d}\\(p-a,b)=1}}1+O(\log x).\]
By the Brun--Titchmarsh inequality, we have
\begin{align*}
\sum_{\substack{d<\sqrt{x_a}\\(d,a)=1}}\frac{\mu(d)^2}{\varphi(d)}\sum_{\substack{p\le d^2+a\\p\equiv a\psmod{d}\\(p-a,b)=1}}&\ll\sum_{d<\sqrt{x_a}}\frac{\mu(d)^2}{\varphi(d)}\cdot\frac{d^2}{\varphi(d)\log 2d}\\
&=\sum_{d\le\sqrt[3]{x_a}}\frac{\mu(d)^2d^2}{\varphi(d)^2\log 2d}+\sum_{\sqrt[3]{x_a}<d<\sqrt{x_a}}\frac{\mu(d)^2d^2}{\varphi(d)^2\log 2d}\\
&\ll\sum_{d\le\sqrt[3]{x_a}}\frac{\mu(d)^2d^2}{\varphi(d)^2}+\frac{1}{\log x_a}\sum_{\sqrt[3]{x_a}<d<\sqrt{x_a}}\frac{\mu(d)^2d^2}{\varphi(d)^2}\\
&\ll\frac{\sqrt{x}}{\log x}.
\end{align*}
Hence, we have
\begin{align}\label{eq:S_102}
S_1(x;a,b)&=\sum_{\substack{d<\sqrt{x_a}\\(d,a)=1}}\frac{\mu(d)^2}{\varphi(d)}\sum_{\substack{p\le x\\p\equiv a\psmod{d}\\(p-a,b)=1}}1+O\left(\frac{\sqrt{x}}{\log x}\right)\nonumber\\
&=\sum_{\substack{c\mid b\\(c,a)=1}}\mu(c)\sum_{\substack{d<\sqrt{x_a}\\(d,a)=1}}\frac{\mu(d)^2}{\varphi(d)}\pi(x;[c,d],a)+O\left(\frac{\sqrt{x}}{\log x}\right)\nonumber\\
&=\sum_{\substack{n<|b|\sqrt{x_a}\\(n,a)=1}}h_n(x;a,b)\pi(x;n,a)+O\left(\frac{\sqrt{x}}{\log x}\right),
\end{align}
where 
\[h_n(x;a,b)\colonequals\sum_{\substack{n=[c,d]\\c\mid b\\d<\sqrt{x_a}\\(c,a)=(d,a)=1}}\frac{\mu(c)\mu(d)^2}{\varphi(d)}.\]
Given positive integers $n<|b|\sqrt{x_a}$, $c\mid b$ with $(c,a)=1$, and $c'\mid c$, there is at most one positive integer $d<\sqrt{x_a}$ with $(d,a)=1$, $(c,d)=c'$ and $[c,d]=n$. If such a $d$ exists, then we must also have $d\asymp n$. Consequently, we have $h_n(x;a,b)\ll 1/\varphi(n)$ for all $n<|b|\sqrt{x_a}$. The Bombieri--Vinogradov theorem \cite[Theorem, \S 28]{Dav} implies that for any fixed $A>0$, 
\begin{equation}\label{eq:BV}
\sum_{d\le Q}\max\limits_{0\le y\le x}\max\limits_{(a,d)=1}\left|\pi(y;d,a)-\frac{\li(y)}{\varphi(d)}\right|\ll \sqrt{x}Q(\log x)^4
\end{equation}
holds uniformly for all $x\ge3$ and all $\sqrt{x}/(\log x)^{A}\le Q\le \sqrt{x}$. Applying \eqref{eq:BV} with $A=A_0+4$ and $Q=\sqrt{x_a}/(\log x_a)^{A_0+4}$, we obtain
\begin{equation}\label{eq:BVh}
\sum_{\substack{n\le\sqrt{x_a}/(\log x_a)^{A_0+4}\\(n,a)=1}}|h_n(x;a,b)|\left|\pi(x;n,a)-\frac{\li(x)}{\varphi(n)}\right|\ll\frac{x}{(\log x)^{A_0}}.
\end{equation}
Since Brun--Titchmarsh implies that
\begin{align*}
\sum_{\substack{\sqrt{x_a}/(\log x_a)^{A_0+4}<n<|b|\sqrt{x_a}\\(n,a)=1}}|h_n(x;a,b)|\left(\frac{\li(x)}{\varphi(n)}+\pi(x;n,a)\right)&\ll\frac{x}{\log x}\sum_{n>\sqrt{x_a}/(\log x_a)^{A_0+4}}\frac{1}{\varphi(n)^2}\\
&\ll\sqrt{x}(\log x)^{A_0+3},
\end{align*}
it follows from \eqref{eq:S_102} and \eqref{eq:BVh} that
\[S_1(x;a,b)=\li(x)\sum_{\substack{n<|b|\sqrt{x_a}\\(n,a)=1}}\frac{h_n(x;a,b)}{\varphi(n)}+O\left(\frac{x}{(\log x)^{A_0}}\right).\]
But
\begin{align*}
\sum_{\substack{n<|b|\sqrt{x_a}\\(n,a)=1}}\frac{h_n(x;a,b)}{\varphi(n)}&=\sum_{\substack{c\mid b\\(c,a)=1}}\mu(c)\sum_{\substack{d<\sqrt{x_a}\\(d,a)=1}}\frac{\mu(d)^2}{\varphi(d)\varphi([c,d])}\\
&=\sum_{\substack{c\mid b\\(c,a)=1}}\mu(c)\sum_{\substack{d\ge1\\(d,a)=1}}\frac{\mu(d)^2}{\varphi(d)\varphi([c,d])}+O\left(\sum_{d\ge\sqrt{x_a}}\frac{\mu(d)^2}{\varphi(d)^2}\right)\\
&=\sum_{\substack{c\mid b\\(c,a)=1}}\mu(c)\sum_{\substack{d\ge1\\(d,a)=1}}\frac{\mu(d)^2}{\varphi(d)\varphi([c,d])}+O\left(\frac{1}{\sqrt{x}}\right),
\end{align*}
whence
\[S_1(x;a,b)=\li(x)\sum_{\substack{c\mid b\\(c,a)=1}}\mu(c)\sum_{\substack{d\ge1\\(d,a)=1}}\frac{\mu(d)^2}{\varphi(d)\varphi([c,d])}+O\left(\frac{x}{(\log x)^{A_0}}\right).\]
To complete the proof of \eqref{eq:S_101}, it suffices to show that the value of the double sum above is equal to $C_{a,b}$. For simplicity of notation, we denote by $D_a$ the first product in the definition \eqref{eq:C_ab} of $C_{a,b}$, i.e., 
\[D_a\colonequals\prod_{p\nmid a}\left(1+\frac{1}{(p-1)^2}\right).\]
Using the identity 
\[\varphi([c,d])=\frac{\varphi(c)\varphi(d)}{\varphi((c,d))},\]
we obtain
\[\sum_{\substack{c\mid b\\(c,a)=1}}\mu(c)\sum_{\substack{d\ge1\\(d,a)=1}}\frac{\mu(d)^2}{\varphi(d)\varphi([c,d])}=\sum_{\substack{c\mid b\\(c,a)=1}}\frac{\mu(c)}{\varphi(c)}\sum_{\substack{d\ge1\\(d,a)=1}}\frac{\mu(d)^2\varphi((c,d))}{\varphi(d)^2}.\]
If $n\in\N$ is square-free, then we may write
\[\varphi(n)=\sum_{m\mid n}\varphi(m)f(m),\]
where
\[f(m)\colonequals\prod_{p\mid m}\left(1-\frac{1}{p-1}\right).\]
Applying this identity with $n=(c,d)$ for square-free $c$ and $d$, we get
\begin{align*}
\sum_{\substack{c\mid b\\(c,a)=1}}\frac{\mu(c)}{\varphi(c)}\sum_{\substack{d\ge1\\(d,a)=1}}\frac{\mu(d)^2}{\varphi(d)\varphi([c,d])}&=\sum_{\substack{c\mid b\\(c,a)=1}}\frac{\mu(c)}{\varphi(c)}\sum_{m\mid c}\varphi(m)f(m)\sum_{\substack{d\ge1\\(d,a)=1\\m\mid d}}\frac{\mu(d)^2}{\varphi(d)^2}\\
&=\sum_{\substack{c\mid b\\(c,a)=1}}\frac{\mu(c)}{\varphi(c)}\sum_{m\mid c}\frac{f(m)}{\varphi(m)}\sum_{\substack{d\ge1\\(d,ma)=1}}\frac{\mu(d)^2}{\varphi(d)^2}\\
&=D_a\sum_{\substack{c\mid b\\(c,a)=1}}\frac{\mu(c)}{\varphi(c)}\sum_{m\mid c}\frac{f(m)}{\varphi(m)}\prod_{\substack{p\mid m\\p\nmid a}}\left(1+\frac{1}{(p-1)^2}\right)^{-1}.
\end{align*}
Let
\[g(n;a)\colonequals \sum_{m\mid n}\frac{f(m)}{\varphi(m)}\prod_{\substack{p\mid m\\p\nmid a}}\left(1+\frac{1}{(p-1)^2}\right)^{-1}.\]
Then $g$ is multiplicative in $n$ and satisfies
\[g(p;a)=1+\frac{f(p)}{\varphi(p)}\left(1+\frac{1}{(p-1)^2}\right)^{-1}=\frac{p(p-1)}{(p-1)^2+1}\]
for all primes $p\nmid a$. It follows that
\begin{align*}
\sum_{\substack{c\mid b\\(c,a)=1}}\frac{\mu(c)}{\varphi(c)}\sum_{\substack{d\ge1\\(d,a)=1}}\frac{\mu(d)^2}{\varphi(d)\varphi([c,d])}&=D_a\sum_{\substack{c\mid b\\(c,a)=1}}\frac{\mu(c)}{\varphi(c)}g(c;a)\\
&=D_a\prod_{\substack{p\mid b\\p\nmid a}}\left(1-\frac{g(p;a)}{p-1}\right)\\
&=D_a\prod_{\substack{p\mid b\\p\nmid a}}\frac{(p-1)(p-2)}{(p-1)^2+1}=C_{a,b}.
\end{align*}
This completes the proof.
\end{proof}

Recall that for any $a\in\Z\setminus\{0\}$, we define 
\[\omega_a^*(n)\colonequals\sum_{\substack{(p-a)\mid n\\p>a}}1.\]
We are now ready to prove the following theorem concerning the first moment of $\omega_a^*$ over shifted primes $p-b$ for any fixed $a,b\in\Z\setminus\{0\}$.

\begin{thm}\label{thm:omega_b^*(p-a)}
For any fixed $a,b\in\Z\setminus\{0\}$, we have
\[\sum_{b<p\le x}\omega_a^*(p-b)=C_{a,b}\frac{x\log\log x}{\log x}+O\left(\frac{x}{\log x}\right)\]
for all $x\ge\max(a,b,3)$, where $C_{a,b}$ is as defined in \eqref{eq:C_ab}.
\end{thm}
\begin{proof}
It is sufficient to prove the theorem for sufficiently large $x$. We have
\begin{equation}\label{eq:omega_a^*(p-b)}
\sum_{b<p\le x}\omega_a^*(p-b)=\sum_{b<p\le x}\sum_{\substack{q-a\mid p-b\\0<q-a<\sqrt{p-b}}}1+\sum_{b<p\le x}\sum_{\substack{n\mid p-b\\n\ge\sqrt{p-b}}}1_{\PP}(n+a).
\end{equation}
The first double sum can be estimated by using Brun--Tichmarsh and \eqref{eq:BV}. Note that
\[\sum_{b<p\le x}\sum_{\substack{q-a\mid p-b\\0<q-a<\sqrt{p-b}}}1=\sum_{\substack{a<q<a+\sqrt{x_b}\\(q-a,b)=1}}\left(\pi(x;q-a,b)-\pi\left((q-1)^2+b;q-a,b\right)\right)+O(\sqrt{x}),\]
where we have once again adopted the shorthand notation $x_b\colonequals x-b$. By Brun--Titchmarsh and Proposition \ref{prop:shiftphi}, we have 
\begin{align*}
\sum_{\substack{a<q<a+\sqrt{x_b}\\(q-a,b)=1}}\pi\left((q-a)^2+a;q-1,b\right)&\ll\sum_{a<q<a+\sqrt{x_b}}\frac{(q-a)^2}{\varphi(q-a)\log q}\\
&\ll\frac{\sqrt{x}}{\log x}\sum_{a<q<a+\sqrt{x_b}}\frac{q-a}{\varphi(q-a)}\\
&\ll\frac{x}{(\log x)^2}.
\end{align*}
It follows that 
\begin{equation}\label{eq:omega_b^*(p-a)sum1}
\sum_{b<p\le x}\sum_{\substack{q-a\mid p-b\\0<q-a<\sqrt{p-b}}}1=\sum_{\substack{a<q<a+\sqrt{x_b}\\(q-a,b)=1}}\pi(x;q-a,b)+O\left(\frac{x}{(\log x)^2}\right).
\end{equation}
Taking $A=5$ and $Q=a+\sqrt{x_b}/(\log x_b)^{5}$ in \eqref{eq:BV}, we see that
\[\sum_{\substack{a<q<a+\sqrt{x_b}/(\log x_b)^5\\(q-a,b)=1}}\pi(x;q-a,b)=\li(x)\sum_{\substack{a<q<a+\sqrt{x_b}/(\log x_b)^5\\(q-a,b)=1}}\frac{1}{\varphi(q-a)}+O\left(\frac{x}{\log x}\right).\]
By Proposition \ref{prop:shiftphi} and partial summation, we obtain
\begin{equation}\label{eq:shift1/phi}
\sum_{\substack{a<q<y\\(q-a,b)=1}}\frac{1}{\varphi(q-a)}=C_{a,b}\log\log y+ O(1)
\end{equation}
for all $y\ge\max(a,3)$. Hence, we have
\[\sum_{\substack{a<q<a+\sqrt{x_b}/(\log x_b)^5\\(q-a,b)=1}}\pi(x;q-a,b)=C_{a,b}\frac{x\log\log x}{\log x}+O\left(\frac{x}{\log x}\right).\]
Since Brun--Titchmarsh and \eqref{eq:shift1/phi} implies that
\[\sum_{\substack{a+\sqrt{x_b}/(\log x_b)^5<q<a+\sqrt{x_b}\\(q-a,b)=1}}\pi(x;q-a,b)\ll\frac{x}{\log x}\sum_{\substack{a+\sqrt{x_b}/(\log x_b)^5<q<a+\sqrt{x_b}\\(q-a,b)=1}}\frac{1}{\varphi(q-a)}\ll\frac{x}{\log x},\]
we conclude that
\[\sum_{\substack{a<q<a+\sqrt{x_b}\\(q-a,b)=1}}\pi(x;q-a,b)=C_{a,b}\frac{x\log\log x}{\log x}+O\left(\frac{x}{\log x}\right).\]
Inserting this estimate into \eqref{eq:omega_b^*(p-a)sum1} yields
\begin{equation}\label{eq:omega_b^*(p-a)sum1.1}
\sum_{b<p\le x}\sum_{\substack{q-a\mid p-b\\0<q-a<\sqrt{p-b}}}1=C_{a,b}\frac{x\log\log x}{\log x}+O\left(\frac{x}{\log x}\right).
\end{equation}
\par To estimate the second double sum in \eqref{eq:omega_a^*(p-b)}, we write $p-b=nd$ and observe that
\begin{equation}\label{eq1:n+a,dn+b}
\sum_{b<p\le x}\sum_{\substack{n\mid p-b\\n\ge\sqrt{p-b}}}1_{\PP}(n+a)\le\sum_{d\le\sqrt{x_b}}\sum_{n\le x_b/d}1_{\PP}(n+a)1_{\PP}(dn+b).
\end{equation}
We may restrict our attention to those $d\le\sqrt{x_b}$ with $\gcd(d,b)=1$, since the inner sum in \eqref{eq1:n+a,dn+b} is clearly at most 1 if $\gcd(d,b)>1$, so that the contribution to \eqref{eq1:n+a,dn+b} from the $d$'s with $\gcd(d,b)>1$ is $\ll \sqrt{x}$. If $d=b/a$, then 
\[\sum_{n\le x_b/d}1_{\PP}(n+a)1_{\PP}(dn+b)\ll 1_{a=b}\cdot\frac{x}{\log x}+1.\]
So the contribution to \eqref{eq1:n+a,dn+b} from $d=b/a$ is $\ll x/\log x$. Suppose now that $d\ne b/a$. By Brun’s or Selberg’s sieve, we have
\[\sum_{n\le x_b/d}1_{\PP}(n+a)1_{\PP}(dn+b)\ll\frac{x}{d(\log x)^2}\prod_{p\mid d(ad-b)}\left(1-\frac{1}{p}\right)^{-1}=\frac{1}{\varphi(d)}\cdot\frac{|ad-b|}{\varphi(|ad-b|)}\cdot\frac{x}{(\log x)^2},\]
since the assumption that $\gcd(d,b)=1$ implies that $\gcd(d,ad-b)=1$. Summing this over all $d\in\N\setminus\{b/a\}$ yields the contribution 
\begin{equation}\label{eq:sbd}
\frac{x}{(\log x)^2}\sum_{\substack{d\le\sqrt{x_b}\\d\ne b/a}}\frac{1}{\varphi(d)}\cdot\frac{|ad-b|}{\varphi(|ad-b|)}.
\end{equation}
By the Cauchy--Schwarz inequality, the sum above is
\[\le\left(\sum_{d\le\sqrt{x_b}}\frac{d}{\varphi(d)^2}\right)^{\frac{1}{2}}\left(\sum_{\substack{d\le\sqrt{x_b}\\d\ne b/a}}\frac{1}{d}\cdot\frac{(ad-b)^2}{\varphi(|ad-b|)^2}\right)^{\frac{1}{2}}.\]
The first factor is easily seen to be $\ll\sqrt{\log x}$. To estimate the second factor, we put $n=ad-b$, so that
\[\sum_{\substack{d\le\sqrt{x_b}\\d\ne b/a}}\frac{1}{d}\cdot\frac{(ad-b)^2}{\varphi(|ad-b|)^2}\ll\sum_{n\le |b|+|a|\sqrt{x_b}}\frac{n}{\varphi(n)^2}\ll\log x.\]
It follows that \eqref{eq:sbd} is $\ll x/\log x$. Combining the above estimates with \eqref{eq1:n+a,dn+b}, we conclude
\begin{equation}\label{eq2:eq:n+a,dn+b}
\sum_{b<p\le x}\sum_{\substack{n\mid p-b\\ n\ge\sqrt{p-b}}}1_{\PP}(n+a)\le\sum_{d\le\sqrt{x_b}}\sum_{n\le x_b/d}1_{\PP}(n+a)1_{\PP}(dn+b)\ll\frac{x}{\log x}.
\end{equation}
Inserting this bound and \eqref{eq:omega_b^*(p-a)sum1.1} into \eqref{eq:omega_a^*(p-b)} completes the proof of the theorem.
\end{proof}
\par It may be of interest to compare the asymptotic formula for the first moment of $\omega^*_a$ over shifted primes $p-b$ supplied by Theorem \ref{thm:omega_b^*(p-a)} with those for $\omega$ and $\tau$, the latter of which resolves the well-known Titchmarsh divisor problem. Assuming GRH, Titchmarsh \cite{Tit} showed that
\begin{equation}\label{TitD}
\sum_{b<p\le x}\tau(p-b)=x\frac{\varphi(|b|)}{|b|}\prod_{p\nmid b}\left(1+\frac{1}{p(p-1)}\right)+O\left(\frac{x\log\log x}{\log x}\right)
\end{equation}
for every fixed $b\in\Z\setminus\{0\}$. Using his dispersion method, Linnik \cite{Lini} provided the first unconditional proof of this asymptotic formula. Soon after the Bombieri--Vinogradov theorem became available, Rodriguez \cite{Rod} and Halberstam \cite{Hal2} obtained independently quick proofs of \eqref{TitD}. In fact, our proof of Theorem \ref{thm:omega_b^*(p-a)} leverages some of the ideas in \cite{Hal2}, and these ideas also lead to the asymptotic formula
\[\sum_{b<p\le x}\omega(p-b)=\frac{x\log\log x}{\log x}+O\left(\frac{x}{\log x}\right)\]
for any given $b\in\Z\setminus\{0\}$. Thus, when averaging over shifted primes, $\omega^*_a$ behaves more like $\omega$ than $\tau$. The Titchmarsh divisor problem and various analogues of it have since been explored extensively. The interested reader may refer to \cite{ABL} for some recent generalizations of \eqref{TitD} and further discussions on this topic, including a uniform version of \eqref{TitD} and an automorphic analogue of the Titchmarsh divisor problem for the Hecke eigenvalues of cuspidal newforms.
\par Given $a,b\in\Z\setminus\{0\}$ such that $2\nmid a$ and $2\mid b$ do not occur simultaneously, Theorem \ref{thm:omega_b^*(p-a)} implies that $\omega_a^*$ is about $C_{a,b}\log\log p$ on average over the shifted primes $p-b$. One may compare this result with the easy fact that $\omega_a^*$ is about $\log\log n$ on average over the positive integers $n$. This analogy becomes more illuminating when we examine the second moments of $\omega_a^*$ over the positive integers and over the shifted primes, respectively. The arguments in \cite{Din,MM} show that
\begin{equation}\label{eq:2momega_a^*}
\frac{1}{x}\sum_{n\le x}\omega_a^*(n)^2=\sum_{\substack{[p-a,q-a]\le x\\p,q>a}}\frac{1}{[p-a,q-a]}+O(1)\asymp\log x
\end{equation}
for any fixed $a\in\Z\setminus\{0\}$, see \cite{FP} also. We shall prove an analogue of this result concerning the second moment of $\omega_a^*$ over shifted primes $p-b$, where $2\mid a$ or $2\nmid b$ are fixed integers. Here the divisibility constraint on $a,b$ is only necessary for obtaining a nontrivial lower bound. 
\par We start with the following (stronger) variant of a theorem of Erd\H{o}s and Prachar \cite{EP}. A generalization of this result will be given in Section \ref{S:f([p-a,q-b])}.

\begin{prop}\label{prop:philcm}
For any fixed $a\in\Z\setminus\{0\}$, we have
\[\sum_{\substack{[p-a,q-a]\le x\\p,q>a}}\frac{[p-a,q-a]}{\varphi([p-a,q-a])}=O(x).\]
\end{prop}
\begin{proof}
Let us write $p-a=uw$ and $q-a=vw$ with $w=\gcd(p-a,q-a)$. Then $\gcd(u,v)=1$ and $[p-a,q-a]=uvw$. Using the inequality $\varphi(mn)\ge\varphi(m)\varphi(n)$ for all $m,n\in\N$, we have
\begin{align}\label{ineq:philcm}
\sum_{\substack{[p-a,q-a]\le x\\p,q>a}}\frac{[p-a,q-a]}{\varphi([p-a,q-a])}&\le\sum_{\substack{uvw\le x\\ (u,v)=1}}\frac{uvw}{\varphi(uv)\varphi(w)}1_{\PP}(uw+a)1_{\PP}(vw+a)\nonumber\\
&=\sum_{\substack{uvw\le x\\ (u,v)=(uvw,a)=1}}\frac{uvw}{\varphi(uv)\varphi(w)}1_{\PP}(uw+a)1_{\PP}(vw+a)+O(x)\nonumber\\
&=\sum_{\substack{uvw\le x\\ u\ne v\\(u,v)=(uvw,a)=1}}\frac{uvw}{\varphi(uv)\varphi(w)}1_{\PP}(uw+a)1_{\PP}(vw+a)+O(x),
\end{align}
where we observe that the contribution from the terms with $u=v$ is clearly $O(x/\log x)$ by Proposition \ref{prop:shiftphi}.
\par Consider first the contribution to \eqref{ineq:philcm} from the terms with the additional constraint $uv\le\sqrt{x}$. By the inequality $n/\varphi(n)\ll \sigma(n)/n$, we have
\begin{align*}
\sum_{w\le x/uv}\frac{w}{\varphi(w)}1_{\PP}(uw+a)1_{\PP}(vw+a)&\ll\sum_{w\le x/uv}\frac{\sigma(w)}{w}1_{\PP}(uw+a)1_{\PP}(vw+a)\\
&\le 2\sum_{w\le x/uv}1_{\PP}(uw+a)1_{\PP}(vw+a)\sum_{\substack{d\mid w\\ d\le \sqrt{w}}}\frac{1}{d}\\
&=2\sum_{d\le \sqrt{x/uv}}\frac{1}{d}\sum_{\substack{w\le x/uv\\ d\mid w}}1_{\PP}(uw+a)1_{\PP}(vw+a)\\
&=2\sum_{d\le \sqrt{x/uv}}\frac{1}{d}\sum_{m\le x/uvd}1_{\PP}(udm+a)1_{\PP}(vdm+a).
\end{align*}
Since $uvd\le \sqrt{xuv}\le x^{3/4}$, Brun's or Selberg's sieve yields
\[\sum_{m\le x/uvd}1_{\PP}(udm+a)1_{\PP}(vdm+a)\ll\frac{x}{uvd(\log x)^2}\prod_{p\mid uvd}\left(1-\frac{1}{p}\right)^{-2}\prod_{\substack{p\mid u-v\\ p\nmid uvd}}\left(1-\frac{1}{p}\right)^{-1}.\]
It follows that
\[\sum_{w\le x/uv}\frac{w}{\varphi(w)}1_{\PP}(uw+a)1_{\PP}(vw+a)\ll\left(\frac{uv}{\varphi(uv)}\right)^2\frac{|u-v|}{\varphi(|u-v|)}\cdot\frac{x}{uv(\log x)^2}.\]
Hence, the contribution to \eqref{ineq:philcm} from the terms with $uv\le\sqrt{x}$ is
\begin{align*}
&\ll\frac{x}{(\log x)^2}\sum_{\substack{uv\le \sqrt{x}\\u\ne v}}\frac{1}{uv}\left(\frac{uv}{\varphi(uv)}\right)^3\frac{|u-v|}{\varphi(|u-v|)}\\
&\le\frac{x}{2(\log x)^2}\sum_{\substack{uv\le \sqrt{x}\\u\ne v}}\frac{1}{uv}\left(\left(\frac{uv}{\varphi(uv)}\right)^6+\left(\frac{u-v}{\varphi(|u-v|)}\right)^2\right).
\end{align*}
It is easy to see that 
\[\sum_{\substack{uv\le \sqrt{x}\\u\ne v}}\frac{1}{uv}\left(\frac{uv}{\varphi(uv)}\right)^6\le\sum_{n\le \sqrt{x}}\frac{\tau(n)}{n}\left(\frac{n}{\varphi(n)}\right)^2\ll(\log x)^2.\]
Moreover, we have, by putting $n=v-u$, that
\begin{align*}
\sum_{\substack{uv\le \sqrt{x}\\u\ne v}}\frac{1}{uv}\left(\frac{u-v}{\varphi(|u-v|)}\right)^2&\le2\sum_{u<v\le\sqrt{x}}\frac{1}{uv}\left(\frac{v-u}{\varphi(v-u)}\right)^2\\
&\le2\sum_{u\le\sqrt{x}}\frac{1}{u}\sum_{n\le\sqrt{x}}\frac{1}{n}\left(\frac{n}{\varphi(n)}\right)^2\ll(\log x)^2.
\end{align*}
Therefore, the contribution to \eqref{ineq:philcm} from the terms with $uv\le\sqrt{x}$ is $O(x)$.
\par Now we estimate the contribution to \eqref{ineq:philcm} from the terms with $uv>\sqrt{x}$. In this case, we must have $w<\sqrt{x}$. Without loss of generality, we may also assume that $u<v$, so that $u<\sqrt{x/w}$ and $uw<\sqrt{xw}<x^{3/4}$. By Brun--Titchmarsh, we have
\begin{align*}
\sum_{\substack{v\le x/uw\\ (v,a)=1}}\frac{v}{\varphi(v)}1_{\PP}(vw+a)&\ll\sum_{\substack{d\le \sqrt{x/uw}\\(d,a)=1}}\frac{1}{d}\,\pi(x/u;wd,a)\\
&\ll\frac{x}{u\log(2x/uw)}\sum_{d\le \sqrt{x/uw}}\frac{1}{d\varphi(wd)}\\
&\ll\frac{x}{u\varphi(w)\log x}.
\end{align*}
It follows that the contribution to \eqref{ineq:philcm} from the terms with $uv>\sqrt{x}$ is 
\[\ll\frac{x}{\log x}\sum_{\substack{w<\sqrt{x}\\ (w,a)=1}}\frac{w}{\varphi(w)^2}\sum_{u<\sqrt{x/w}}\frac{1}{\varphi(u)}1_{\PP}(uw+a).\]
Recall that $p=uw+a$, so that the expression above is 
\begin{align*}
&\le\frac{x}{\log x}\sum_{w<\sqrt{x}}\frac{w}{\varphi(w)^2}\sum_{\substack{a<p\le a+\sqrt{xw}\\ p\equiv a\psmod{w}}}\frac{1}{\varphi((p-a)/w)}\\
&=\frac{x}{\log x}\sum_{w<\sqrt{x}}\frac{w^2}{\varphi(w)^2}\sum_{\substack{a<p\le a+\sqrt{xw}\\ p\equiv a\psmod{w}}}\frac{1}{p-a}\cdot\frac{(p-a)/w}{\varphi((p-a)/w)}\\
&\le\frac{x}{\log x}\sum_{w<\sqrt{x}}\frac{w^2}{\varphi(w)^2}\sum_{\substack{a<p\le a+\sqrt{xw}\\ p\equiv a\psmod{w}}}\frac{1}{\varphi(p-a)}\\
&\le\frac{x}{\log x}\sum_{a<p\le a+x}\frac{h(p-a)}{\varphi(p-a)},
\end{align*}
where
\[h(n)\colonequals\sum_{d\mid n}\frac{d^2}{\varphi(d)^2}\]
for every $n\in\N$, and the third inequality is an immediate consequence of the fact that the function $n/\varphi(n)$ is increasing on the divisor lattice on $\N$. By \cite[Theorem 1]{Poll}, we have
\[\sum_{a<p\le a+x^{3/4}}h(p-a)\frac{p-a}{\varphi(p-a)}\ll x,\]
from which it follows that
\[\sum_{a<p\le a+x}\frac{h(p-a)}{\varphi(p-a)}\ll\log x.\]
Hence, we conclude that the contribution to \eqref{ineq:philcm} from the terms with $uv>\sqrt{x}$ is $O(x)$. Combining this with the contribution to \eqref{ineq:philcm} from the terms with $uv\le\sqrt{x}$ finishes the proof of Proposition \ref{prop:philcm}.
\end{proof}

We are now in a position to prove the following analogue of \eqref{eq:2momega_a^*} on the second moment of $\omega_a^*$ on shifted primes $p-b$. Notably, we shall not only apply Proposition \ref{prop:philcm} to handle the upper bound but also recycle the ideas used in its proof presented above. An easy variant of \eqref{eq:BV} that we shall use is 
\begin{equation}\label{eq:BVtau}
\sum_{d\le Q}\tau_{\kappa}(d)\max\limits_{0\le y\le x}\max\limits_{(a,d)=1}\left|\pi(y;d,a)-\frac{\li(y)}{\varphi(d)}\right|\ll x^{3/4}\sqrt{Q}(\log x)^{(\kappa^2+3)/2}
\end{equation}
uniformly for all $x\ge3$ and all $\sqrt{x}/(\log x)^{A}\le Q\le \sqrt{x}$, where $\kappa,A>0$ are fixed. It is not hard to see that \eqref{eq:BVtau} follows readily from \eqref{eq:BV}. Indeed, if we write
\[\Delta(y;d,a)\colonequals\pi(y;d,a)-\frac{\li(y)}{\varphi(d)},\]
then Cauchy--Schwarz implies that
\[\sum_{d\le Q}\tau_{\kappa}(d)\max\limits_{0\le y\le x}\max\limits_{(a,d)=1}|\Delta(y;d,a)|\le\left(\sum_{d\le Q}\frac{\tau_{\kappa}(d)^2}{\varphi(d)}\right)^{\frac{1}{2}}\left(\sum_{d\le Q}\varphi(d)\max\limits_{0\le y\le x}\max\limits_{(a,d)=1}\Delta(y;d,a)^2\right)^{\frac{1}{2}}.\]
The first factor on the right-hand side is $\ll(\log x)^{\kappa^2/2}$, whereas the second factor on the right-hand side is 
\[\ll\left(\frac{x}{\log x}\sum_{d\le Q}\max\limits_{0\le y\le x}\max\limits_{(a,d)=1}|\Delta(y;d,a)|\right)^{\frac{1}{2}}\ll x^{3/4}\sqrt{Q}(\log x)^{3/2}\]
by \eqref{eq:BV}, since Brun--Titchmarsh implies that
\[\Delta(y;d,a)\ll1_{>y}(d)+1_{\le y}(d)\cdot\frac{y}{\varphi(d)\log(2y/d)}\ll\frac{x}{\varphi(d)\log x}\]
for $y\in[0,x]$ and $d\le Q$. Thus \eqref{eq:BVtau} follows.

\begin{thm}\label{thm:2mshiftomega_a^*}
For any fixed $a,b\in\Z\setminus\{0\}$ such that $2\mid a$ or $2\nmid b$, we have
\[\frac{1}{\pi(x)}\sum_{b<p\le x}\omega_a^*(p-b)^2\asymp\log x\]
for all $x\ge\max(a,b,3)$.
\end{thm}
\begin{proof}
We may assume that $x$ is sufficiently large. The starting point of the proof of the asserted lower bound is the simple observation that
\begin{align}\label{eq:2mshiftomega_a^*lb}
\sum_{b<r\le x}\omega_a^*(r-b)^2&=\sum_{\substack{[p-a,q-a]\le x\\p,q>a}}(\pi(x;[p-a,q-a],b)-\pi(b;[p-a,q-a],b))\nonumber\\
&\ge\sum_{\substack{[p-a,q-a]\le x^{1/3}\\p,q>a}}\pi(x;[p-a,q-a],b)+O\left(x^{1/3}\right),
\end{align}
where we have used the fact that 
\begin{equation}\label{eq:EPa}
\#\left\{(p,q)\in(\N\cap(a,\infty))^2\colon[p-a,q-a]\le x\right\}=O(x),
\end{equation}
which is a direct consequence of Proposition \ref{prop:philcm}. Let 
\[\beta_a(n)\colonequals\#\left\{(p,q)\in(\N\cap(a,\infty))^2\colon[p-a,q-a]=n\right\}\]
for every $n\in\N$. Then we have 
\begin{equation}\label{eq:upi1}
\sum_{\substack{[p-a,q-a]\le x^{1/3}\\p,q>a}}\pi(x;[p-a,q-a],b)=\sum_{n\le x^{1/3}}\beta_a(n)\pi(x;n,b).
\end{equation}
Since $\beta_a(n)\le\tau_3(n)$, it follows that
\begin{align}\label{eq:upi2}
\sum_{n\le x^{1/3}}\beta_a(n)\pi(x;n,b)&=\sum_{\substack{n\le x^{1/3}\\(n,b)=1}}\beta_a(n)\pi(x;n,b)+O\left(\sum_{n\le x^{1/3}}\tau_3(n)\right)\nonumber\\
&=\sum_{\substack{n\le x^{1/3}\\(n,b)=1}}\beta_a(n)\pi(x;n,b)+O\left(x^{1/3}(\log x)^2\right).
\end{align}
To estimate the sum above, we appeal to \eqref{eq:BVtau} with $\kappa=3$, $A=14$, and $Q=\sqrt{x}/(\log x)^{14}$ to obtain

\[\sum_{\substack{n\le x^{1/3}\\(n,b)=1}}\beta_a(n)\left|\pi(x;n,b)-\frac{\li(x)}{\varphi(n)}\right|\ll\frac{x}{\log x}.\]
Hence, we find that
\begin{align*}
\sum_{\substack{n\le x^{1/3}\\(n,b)=1}}\beta_a(n)\pi(x;n,b)&=\li(x)\sum_{\substack{n\le x^{1/3}\\(n,b)=1}}\frac{\beta_a(n)}{\varphi(n)}+O\left(\frac{x}{\log x}\right)\nonumber\\
&\ge\li(x)\sum_{\substack{n\le x^{1/3}\\(n,b)=1}}\frac{\beta_a(n)}{n}+O\left(\frac{x}{\log x}\right)\nonumber\\
&=\li(x)\sum_{\substack{[p-a,q-a]\le x^{1/3}\\p,q>a\\(p-a,b)=(q-a,b)=1}}\frac{1}{[p-a,q-a]}+O\left(\frac{x}{\log x}\right).
\end{align*}
Combining this estimate with \eqref{eq:2mshiftomega_a^*lb}--\eqref{eq:upi2}, we see that the lower bound asserted in the theorem will follow if we can show 
\begin{equation}\label{eq:2momega_ab}
\sum_{\substack{[p-a,q-a]\le x\\p,q>a\\(p-a,b)=(q-a,b)=1}}\frac{1}{[p-a,q-a]}\gg\log x.
\end{equation}
This lower bound, which is slightly more general than the one supplied by \eqref{eq:2momega_a^*}, can be easily verified by adapting the proof of \cite[Eq. (29)]{FP}. Assume that $b$ is square-free. As in the proof of \cite[Eq. (29)]{FP}, we start by observing that 
\begin{align*}
\sum_{\substack{[p-a,q-a]\le x\\p,q>a\\(p-a,b)=(q-a,b)=1}}\frac{1}{[p-a,q-a]}&\gg\sum_{\substack{d\le x^{\delta/4}/|b|\\(d,ab)=1}}\mu(d)^2\varphi(d)\sum_{\substack{x^{1/4}<p,q\le\sqrt{x}+a\\p,q\equiv a\psmod{d}\\(p-a,b)=(q-a,b)=1}}\frac{1}{pq}\\
&\ge\sum_{\substack{c\in(\Z/b\Z)^{\times}\\a+c\in(\Z/b\Z)^{\times}}}\sum_{\substack{d\le x^{\delta/4}/|b|\\(d,ab)=1}}\mu(d)^2\varphi(d)\left(\sum_{\substack{x^{1/4}<p\le\sqrt{x}+a\\p\equiv a\psmod{d}\\p\equiv a+c\psmod{b}}}\frac{1}{p}\right)^2
\end{align*}
uniformly for all sufficiently large $x$ and all $\delta\in(0,1)$. Our assumption that $2\mid a$ or $2\nmid b$ ensures that there exists $c\in(\Z/b\Z)^{\times}$ with $a+c\in(\Z/b\Z)^{\times}$. Indeed, by the Chinese Remainder Theorem, it suffices to prove this when $b=p^{\nu}$ is a prime power. If $p>2$, then any $c\in\Z$ that avoids the residue classes $0$ and $-a$ $\psmod{p}$ works. If $p=2$, then our assumption implies that $2\mid a$, so that we can take $c$ to be any odd integer. This verifies our claim. Now if we denote by $k(a,b,c,d)$ the unique element of $(\Z/bd\Z)^{\times}$ satisfying $k(a,b,c,d)\equiv a\psmod{d}$ and $k(a,b,c,d)\equiv a+c\psmod{b}$, then 
\[\sum_{\substack{[p-a,q-a]\le x\\p,q>a\\(p-a,b)=(q-a,b)=1}}\frac{1}{[p-a,q-a]}\gg\sum_{\substack{c\in(\Z/b\Z)^{\times}\\a+c\in(\Z/b\Z)^{\times}}}\sum_{\substack{d\le x^{\delta/4}/|b|\\(d,ab)=1}}\mu(d)^2\varphi(d)\left(\sum_{\substack{x^{1/4}<p\le\sqrt{x}+a\\p\equiv k(a,b,c,d)\psmod{bd}}}\frac{1}{p}\right)^2.\]
Applying \cite[Corollary 1]{FP} as in the proof of \cite[Eq. (29)]{FP}, one finds that 
\[\sum_{\substack{x^{1/4}<p\le\sqrt{x}+a\\p\equiv k(a,b,c,d)\psmod{bd}}}\frac{1}{p}\gg\frac{1}{\varphi(b)\varphi(d)}\]
for some suitable choice of $\delta$, provided that $bd$ is not divisible by a certain prime $s(x^{1/4})>(1/4)\log\log x$. It follows that 
\[\sum_{\substack{[p-a,q-a]\le x\\p,q>a\\(p-a,b)=(q-a,b)=1}}\frac{1}{[p-a,q-a]}\gg\sum_{\substack{d\le x^{\delta/4}/|b|\\(d,s(x^{1/4})ab)=1}}\frac{\mu(d)^2}{\varphi(d)}\gg\log x,\]
completing the proof of \eqref{eq:2momega_ab}.
\par Now we prove the asserted upper bound. We have
\begin{equation}\label{2mshiftomega_a^*ub}
\sum_{b<r\le x}\omega_a^*(r-b)^2=\sum_{b<r\le x}\sum_{\substack{[p-a,q-a]\mid r-b\\p,q>a\\ [p-a,q-a]\le\sqrt{r-b}}}1+\sum_{b<r\le x}\sum_{\substack{n\mid r-b\\ n>\sqrt{r-b}}}\beta_a(n).
\end{equation}
To estimate the first double sum, we apply Brun--Titchmarsh to get
\begin{align}\label{2mshiftomega_a^*ub1}
\sum_{b<r\le x}\sum_{\substack{[p-a,q-a]\mid r-b\\p,q>a\\ [p-a,q-a]\le\sqrt{r-b}}}1&\le\sum_{\substack{[p-a,q-a]\le \sqrt{x_b}\\(p-a,q-a,b)=1\\p,q>a}}\pi(x;[p-a,q-a],b)+O\left(\sqrt{x}\right)\nonumber\\
&\ll\frac{x}{\log x}\sum_{\substack{[p-a,q-a]\le \sqrt{x_b}\\p,q>a}}\frac{1}{\varphi([p-a,q-a])}+O\left(\sqrt{x}\right)\ll x,
\end{align}
where we have also used \eqref{eq:EPa} in the first inequality and Proposition \ref{prop:philcm} in the last inequality. For the second double sum, we have, by writing $r-b=nd$, that
\begin{align*}
\sum_{b<r\le x}\sum_{\substack{n\mid r-b\\ n>\sqrt{r-b}}}\beta_a(n)&\le\sum_{d<\sqrt{x_b}}\sum_{n\le x_b/d}\beta_a(n)1_{\PP}(dn+b)\\
&=\sum_{\substack{d<\sqrt{x_b}\\(d,b)=1}}\sum_{\substack{n\le x_b/d\\(n,b)=1}}\beta_a(n)1_{\PP}(dn+b)+O\left(x^{1/2+o(1)}\right),
\end{align*}
since $\beta_a(n)\le\tau_3(n)=n^{o(1)}$. Furthermore, Brun--Titchmarsh and the inequality $\beta_a(n)\ll 1$ for all $n\in\N$ with $\gcd(n,a)>1$ imply that 
\[\sum_{\substack{d<\sqrt{x_b}\\(d,b)=1}}\sum_{\substack{n\le x_b/d\\(n,b)=1\\(n,a)>1}}\beta_a(n)1_{\PP}(dn+b)\ll\frac{x}{\log x}\sum_{d<\sqrt{x_b}}\frac{1}{\varphi(d)}\ll x.\]
It follows that 
\[\sum_{b<r\le x}\sum_{\substack{n\mid r-b\\ n>\sqrt{r-b}}}\beta_a(n)\le\sum_{\substack{d<\sqrt{x_b}\\(d,b)=1}}\sum_{\substack{n\le x_b/d\\(n,ab)=1}}\beta_a(n)1_{\PP}(dn+b)+O(x).\]
Recalling the definition
\[\beta_a(n)=\sum_{\substack{uvw=n\\(u,v)=1}}1_{\PP}(uw+a)1_{\PP}(vw+a),\]
we arrive at
\begin{equation}\label{eq:sum2beta_a}
\sum_{b<r\le x}\sum_{\substack{n\mid r-b\\ n>\sqrt{r-b}}}\beta_a(n)\le \sum_{\substack{d<\sqrt{x_b}\\(d,b)=1}}\sum_{\substack{uvw\le x_b/d\\(u,v)=1\\(uvw,ab)=1}}1_{\PP}(uw+a)1_{\PP}(vw+a)1_{\PP}(duvw+b)+O(x).
\end{equation}
There are three boundary cases: $u=v$, $adu=b$, and $adv=b$. In fact, the latter two cases simplify to $d=u=1$ and $d=v=1$, respectively, due to the constraint $\gcd(duv,b)=1$. By \eqref{eq2:eq:n+a,dn+b}, the contribution to \eqref{eq:sum2beta_a} from the boundary case $u=v$ is bounded above by
\begin{equation}\label{eq:sum2beta_abd1}
\sum_{d<\sqrt{x_b}}\sum_{w\le x_b/d}1_{\PP}(w+a)1_{\PP}(dw+b)+O(x)\ll x.
\end{equation}
On the other hand, we see from \eqref{eq:EPa} that the total contribution to \eqref{eq:sum2beta_a} from the boundary cases $adu=b$ and $adv=b$ is at most
\begin{equation}\label{eq:sum2beta_abd2}
2\sum_{vw\le x_b}1_{\PP}(w+a)1_{\PP}(vw+a)+O(x)\le2\sum_{n\le x_b}\beta_{a}(n)+O(x)\ll x.
\end{equation}
\par With the boundary cases handled, we proceed to estimate \eqref{eq:sum2beta_a} in the non-boundary case where $u\ne v$, $adu\ne b$, and $adv\ne b$. In this case, we have
\[D=D(a,b,d,u,v)\colonequals aduv(u-v)(adu-b)(adv-b)\ne0.\]
First of all, we consider the contribution to \eqref{eq:sum2beta_a} from the terms subject to the constraint $uv\le(x_b/d)^{2/3}$, so that $uvd\le x_b^{2/3}d^{1/3}\le x_b^{5/6}$. Brun's or Selberg's sieve shows that this contribution is 
\begin{equation}\label{eq:uv<sqrt{x_b/d}1}
\ll\frac{x}{(\log x)^3}\sum_{\substack{d<\sqrt{x_b}\\(d,b)=1}}\sum_{\substack{uv\le (x_b/d)^{2/3}\\u\ne v,\,(u,v)=1\\(uv,ab)=1\\ D\ne0}}\frac{1}{uvd}\prod_{p\mid D}\left(1-\frac{1}{p}\right)^{-2}.
\end{equation}
Taking advantage of the symmetry between $u$ and $v$, we find that the double sum above is
\begin{align}\label{eq:uv<sqrt{x_b/d}2}
&\ll\sum_{\substack{d<x_b\\(d,b)=1}}\sum_{\substack{uv\le x_b/d\\u\ne v,\,D\ne0}}\frac{1}{uvd}\left(\frac{uvd}{\varphi(uvd)}\right)^2\left(\frac{v-u}{\varphi(v-u)}\right)^2\left(\frac{adu-b}{\varphi(|adu-b|)}\right)^2\left(\frac{adv-b}{\varphi(|adv-b|)}\right)^2\nonumber\\
&\le\frac{1}{4}\sum_{\substack{d<x_b\\(d,b)=1}}\sum_{\substack{uv\le x_b/d\\u\ne v,\,adu\ne b}}\frac{1}{uvd}\left(\left(\frac{uvd}{\varphi(uvd)}\right)^8+\left(\frac{v-u}{\varphi(v-u)}\right)^8+2\left(\frac{adu-b}{\varphi(|adu-b|)}\right)^8\right).
\end{align}
We first note that
\begin{equation}\label{eq:uv<sqrt{x_b/d}1a}
\sum_{\substack{d<x_b\\(d,b)=1}}\sum_{\substack{uv\le x_b/d\\u\ne v,\,adu\ne b}}\frac{1}{uvd}\left(\frac{uvd}{\varphi(uvd)}\right)^8\le\sum_{n\le x_b}\frac{\tau_3(n)}{n}\left(\frac{n}{\varphi(n)}\right)^8\ll(\log x)^3.
\end{equation}
Next, we have, by putting $n=v-u$, that
\begin{align*}
\sum_{\substack{uv\le x_b/d\\u\ne v}}\frac{1}{uv}\left(\frac{v-u}{\varphi(v-u)}\right)^8&=2\sum_{\substack{uv\le x_b/d\\u< v}}\frac{1}{uv}\left(\frac{v-u}{\varphi(v-u)}\right)^8\\
&\le2\sum_{u\le \sqrt{x_b}}\frac{1}{u}\sum_{n\le x_b}\frac{1}{n}\left(\frac{n}{\varphi(n)}\right)^8\ll(\log x)^2,
\end{align*}
which implies that
\begin{equation}\label{eq:uv<sqrt{x_b/d}1b}
\sum_{\substack{d<x_b\\(d,b)=1}}\sum_{\substack{uv\le x_b/d\\u\ne v,\,adu\ne b}}\frac{1}{uvd}\left(\frac{v-u}{\varphi(v-u)}\right)^8\ll(\log x)^2\sum_{d\le x_b}\frac{1}{d}\ll(\log x)^3.
\end{equation}
Finally, a similar argument shows that
\begin{align*}
&\hspace*{5mm}\sum_{\substack{d<x_b\\(d,b)=1}}\sum_{\substack{uv\le x_b/d\\u\ne v,\,adu\ne b}}\frac{1}{uvd}\left(\frac{adu-b}{\varphi(|adu-b|)}\right)^8\\
&\ll\log x\sum_{\substack{d<x_b\\(d,b)=1}}\sum_{\substack{u\le x_b/d\\ adu\ne b}}\frac{1}{ud}\left(\frac{adu-b}{\varphi(|adu-b|)}\right)^8\\
&\le\log x\sum_{\substack{n\le x_b\\ n\ne b/a}}\frac{\tau(n)}{n}\left(\frac{an-b}{\varphi(|an-b|)}\right)^8.
\end{align*}
By \cite[Corollary 3]{NT} with $F_1(n)=\tau(n)$, $F_2(n)=(n/\varphi(n))^8$, $Q_1(t)=t$ and $Q_2(t)=at-b$, we find that 
\[\sum_{\substack{y<n\le 2y\\ n\ne b/a}}\tau(n)\left(\frac{an-b}{\varphi(|an-b|)}\right)^8\ll y\prod_{p\le y}\left(1-\frac{2}{p}\right)\left(\sum_{n\le y}\frac{\tau(n)}{n}\right)\left(\sum_{n\le y}\left(\frac{n}{\varphi(n)}\right)^8\right)\ll y\log y\]
for all $y\ge y_0\colonequals \max(b/a,2)$. Applying this inequality with $y=x_b/2^{j}$ for each $ j\in\N\cap[1,\log(x_b/y_0)/\log 2]$ and summing up the resulting estimates, we obtain
\[\sum_{\substack{n\le x_b\\ n\ne b/a}}\tau(n)\left(\frac{an-b}{\varphi(|an-b|)}\right)^8\ll x\log x.\]
It follows by partial summation that
\begin{equation}\label{eq:uv<sqrt{x_b/d}1c}
\sum_{\substack{d<x_b\\(d,b)=1}}\sum_{\substack{uv\le x_b/d\\u\ne v,\,adu\ne b}}\frac{1}{uvd}\left(\frac{adu-b}{\varphi(|adu-b|)}\right)^8\ll\log x\sum_{\substack{n\le x_b\\ n\ne b/a}}\frac{\tau(n)}{n}\left(\frac{an-b}{\varphi(|an-b|)}\right)^8\ll (\log x)^3.
\end{equation}
Combining \eqref{eq:uv<sqrt{x_b/d}1}--\eqref{eq:uv<sqrt{x_b/d}1c}, we conclude that the contribution to \eqref{eq:sum2beta_a} in the non-boundary case from the terms with the additional constraint $uv\le(x_b/d)^{2/3}$ is $O(x)$.
\par It remains to estimate the contribution to \eqref{eq:sum2beta_a} in the non-boundary case from the terms when $uv>(x_b/d)^{2/3}$. In this case, we must have $w<\sqrt[3]{x_b/d}$. By symmetry, we may also add the constraint that $u<v$, so that $u<\sqrt{x_b/d}$ as well. Again, since $duw\le x_b^{5/6}d^{1/6}\le x_b^{11/12}$, classical results from sieve theory yield
\[\sum_{\substack{v\le x_b/(duw)\\v\ne u,\,(v,u)=1\\adv\ne b,\,(v,ab)=1}}1_{\PP}(vw+a)1_{\PP}(duvw+b)\ll\frac{x}{duw(\log x)^2}\prod_{p\mid w}\left(1-\frac{1}{p}\right)^{-2}\prod_{\substack{p\mid du(dua-b)\\ p\nmid w}}\left(1-\frac{1}{p}\right)^{-1}.\]
Thus, the contribution to \eqref{eq:sum2beta_a} in the non-boundary case from the terms with the constraint $uv>(x_b/d)^{2/3}$ is 
\[\ll\frac{x}{(\log x)^2}\sum_{\substack{d<\sqrt{x_b}\\(d,b)=1}}\sum_{\substack{u<\sqrt{x_b/d}\\adu\ne b,\,(u,ab)=1}}\frac{1}{du}\cdot\frac{du|dua-b|}{\varphi(du|dua-b|)}\sum_{\substack{w<\sqrt[3]{x_b/d}\\(w,ab)=1}}\frac{1}{w}\left(\frac{w}{\varphi(w)}\right)^21_{\PP}(uw+a).\]
Since
\begin{align*}
\sum_{\substack{w<\sqrt[3]{x_b/d}\\(w,ab)=1}}\frac{1}{w}\left(\frac{w}{\varphi(w)}\right)^21_{\PP}(uw+a)&\ll u\sum_{\substack{a+u<p<a+u\sqrt[3]{x_b}\\ p\equiv a\psmod{u}}}\frac{1}{p}\left(\frac{(p-a)/u}{\varphi((p-a)/u)}\right)^2\\
&\le u\sum_{\substack{a<p\le a+x_b\\p\equiv a\psmod{u}}}\frac{1}{p}\left(\frac{p-a}{\varphi(p-a)}\right)^2,
\end{align*}
the sieve bound above is
\[\ll\frac{x}{(\log x)^2}\sum_{a<p\le a+x_b}\frac{\tau(p-a)}{p}\left(\frac{p-a}{\varphi(p-a)}\right)^2\sum_{\substack{n<x_b\\n\ne b/a,\,(n,b)=1}}\frac{1}{n}\cdot\frac{n}{\varphi(n)}\cdot\frac{|an-b|}{\varphi(|an-b|)}.\]
Applying \cite[Corollary 3]{NT} with  $F_1(n)=F_2(n)=n/\varphi(n)$, $Q_1(t)=t$ and $Q_2(t)=at-b$ and arguing as in the proof of \eqref{eq:uv<sqrt{x_b/d}1c}, we have
\[\sum_{\substack{n<x_b\\n\ne b/a,\,(n,b)=1}}\frac{1}{n}\cdot\frac{n}{\varphi(n)}\cdot\frac{|an-b|}{\varphi(|an-b|)}\ll\log x.\]
Hence, the contribution to \eqref{eq:sum2beta_a} in the non-boundary case from the terms with the constraint $uv>(x_b/d)^{2/3}$ is 
\[\ll\frac{x}{\log x}\sum_{a<p\le a+x_b}\frac{\tau(p-a)}{p}\left(\frac{p-a}{\varphi(p-a)}\right)^2\ll x,\]
where the sum above is easily seen to be $O(\log x)$ by \cite[Theorem 1]{Poll} together with partial summation.
\par In conclusion, we have shown that the entire contribution to \eqref{eq:sum2beta_a} in the non-boundary case is $O(x)$. Adding up this contribution and the contributions to \eqref{eq:sum2beta_a} in the boundary cases supplied by \eqref{eq:sum2beta_abd1} and \eqref{eq:sum2beta_abd2} shows that \eqref{eq:sum2beta_a} is $O(x)$. Theorem \ref{thm:2mshiftomega_a^*} follows upon inserting this and \eqref{2mshiftomega_a^*ub1} into \eqref{2mshiftomega_a^*ub}.
\end{proof}

In view of \eqref{eq:2momega_a^*} and Theorem \ref{thm:2mshiftomega_a^*}, it is natural to conjecture that if $a,b\in\Z\setminus\{0\}$ are fixed such that $2\mid a$ or $2\nmid b$, then
\[\frac{1}{\pi(x)}\sum_{b<r\le x}\omega_a^*(r-b)^2\sim \sum_{\substack{[p-a,q-a]\le x\\p,q>a}}\frac{1}{\varphi([p-a,q-a])}\sim c(a,b)\log x\]
for some constant $c(a,b)>0$. More generally, we conjecture that for any fixed $a,b\in\Z\setminus\{0\}$ such that $2\mid a$ or $2\nmid b$ and for any $k\ge2$ , there exists a constant $c_k(a,b)>0$ such that
\[\frac{1}{\pi(x)}\sum_{b<r\le x}\omega_a^*(r-b)^k\sim c_k(a,b)(\log x)^{2^k-k-1}.\]
As mentioned earlier, a similar conjecture on the $k$th moment of $\omega^*$ over natural numbers was made in \cite{FP}. 

\medskip
\section{Mean values of sub-multiplicative functions over shifted primes}\label{S:f([p-a,q-b])}
In this section, we study 
\[\sum_{\substack{[p-a,q-b]\le x\\p>a,q>b}}f([p-a,q-b]),\]
which may be viewed as a two-dimensional analogue of the following sum over shifted primes:
\[\sum_{a<p\le x}f(p-a),\]
where $f$ is a nonnegative-valued arithmetic function, and $a,b\in\Z\setminus\{0\}$ are fixed. As we have seen, this sum is intimately related to the second moments of $\omega_a^*$ over natural numbers and over shifted primes. We shall prove a uniform upper bound valid for a nice class of sub-multiplicative functions, generalizing Proposition \ref{prop:philcm}. Our main tool is a variant of Shiu's theorem due to Pollack \cite{Poll} on the mean values of nonnegative-valued sub-multiplicative functions over a sifited set. We state his result as follows. An arithmetic function $f\colon\N\to\R$ is called {\it sub-multiplicative} if $f(mn)\le f(m)f(n)$ whenever $\gcd(m,n)=1$. Let $A_1>0$ be an absolute constant, and let $A_2\colon\R_{>0}\to\R_{>0}$ be any function. Denote by $\mathscr{M}_s(A_1,A_2)$ the class of sub-multiplicative functions $f\colon \N\to\R_{\ge0}$ satisfying the following two conditions:
\begin{enumerate}[label=(\roman*)]
\item $f(n) \le A_1^{\Omega(n)}$ for all $n\in\N$.
\item Given every $\epsilon > 0$, one has $f(n)\le A_2(\epsilon)n^{\epsilon}$ for all $n\in\N$.
\end{enumerate}
Then \cite[Theorem 1]{Poll} (together with Remark (ii) in \cite{Poll}) asserts that\footnote{Although \cite[Theorem 1]{Poll} is stated for multiplicative functions, its proof shows that sub-multiplicativity suffices. In addition, the theorem is still valid if one assumes only that (ii) holds for a particular value of $\epsilon$ depending on $\alpha$. The same observations have also been made in \cite{NT} on Shiu's theorem. 
} if $\alpha\in(0,1/2)$, $k\in\Z_{\ge0}$, and $f\in\mathscr{M}_s(A_1,A_2)$, and if $\mathscr{N}$ is the subset of $\N\cap(x-y,x]$ of integers whose reductions modulo $p$ avoid a subset $\mathscr{E}_p\subseteq\Z/p\Z$ of cardinality $\#\mathscr{E}_p=\nu(p)\le k$ for every prime $p\le x$, where $x^{\alpha}<y\le x$, then 
\begin{equation}\label{eq:Poll}
\sum_{n\in\mathscr{N}}f(n)\ll_{\alpha, k, A_1,A_2}C_{f,\mathscr{Q}}\frac{y}{\log 2x}\exp\left(\sum_{p\le x}\frac{f(p)-\nu(p)}{p}\right)
\end{equation}
for all $x\ge 1$, where $\mathscr{Q}\colonequals\{A_1<p\le x\colon 0\in\mathscr{E}_p\}$ and 
\[C_{f,\mathscr{Q}}\colonequals\prod_{p\in\mathscr{Q}}\left(1-\frac{f(p)}{p}\right)\left(1-\frac{1}{p}\right)^{-1}.\]
We deduce from \eqref{eq:Poll} the following extension.
\begin{lem}\label{lem:Pollf(mn)}
Let $\alpha,\delta_0\in(0,1)$, $k\in\Z_{\ge0}$, $\eta_0,A_1>0$, $A_2\colon \R_{>0}\to\R_{>0}$, and $f\in\mathscr{M}_s(A_1,A_2)$. If $\mathscr{N}$ is a subset of $\N\cap(x-y,x]$ of integers whose reductions modulo $p$ avoid a subset $\mathscr{E}_p\subseteq\Z/p\Z$ of cardinality $\#\mathscr{E}_p=\nu(p)\le k$ for every prime $p\le x$, where $x^{\alpha}<y\le x$, then
\[\sum_{n\in\mathscr{N}}f(mn)\ll_{\alpha,k,\delta_0,\eta_0,A_1,A_2}C_{f,\mathscr{Q}}(m)m\frac{y}{\log 2x}\exp\left(\sum_{p\le x}\frac{f(p)-\nu(p)}{p}\right)\prod_{p^v\parallel m}\sum_{\nu\ge v}\frac{f(p^{\nu})}{p^{\nu}}+y^{\delta_0}\]
for all $x\ge 1$ and all $m\in\N\cap[1,x^{\eta_0}]$, where
\[C_{f,\mathscr{Q}}(m)\colonequals C_{f,\mathscr{Q}}\prod_{\substack{p>\max(A_1,k)\\p\mid m,p\notin\mathscr{Q}}}\left(1-\frac{f(p)}{p}\right)\left(1-\frac{\nu(p)}{p}\right)^{-1}.\]
\end{lem}
\begin{proof}
It suffices to prove the lemma for sufficiently large $x$ depending on $\alpha,k,\delta_0,\eta_0,A_1,A_2$. In the proof, all the implicit constants will also depend at most on these parameters. Let $\epsilon\colonequals\alpha\delta_0/(10+2\eta_0)\in(0,\alpha/10)$, and put $z\colonequals (y/x^{\epsilon})^{1/(1-\epsilon)}$ and $P_{\mathscr{E}}\colonequals\prod_{p: 0\in\mathscr{E}_p}p$. Then
\[\sum_{n\in\mathscr{N}}f(mn)\le\sum_{\substack{n_1\le x\\ \rad(n_1)\mid m}}f(mn_1)\sum_{\substack{ n_1n_2\in\mathscr{N}\\ (n_2,m)=1}}f(n_2)=\sum_{\substack{n_1\le x\\ \rad(n_1)\mid m\\ (n_1,P_{\mathscr{E}})=1}}f(mn_1)\sum_{\substack{ n_1n_2\in\mathscr{N}\\ (n_2,m)=1}}f(n_2),\]
since the inner sum over $n_2$ vanishes if there exists $p\mid n_1$ such that $0\in\mathscr{E}_p$. If $n_1<z$, then we have $y/n_1>(x/n_1)^{\epsilon}$. By \eqref{eq:Poll}, the inner sum over $n_2$ is
\begin{align}\label{eq:sumn_2}
&\ll\frac{y}{n_1\log(2x/n_1)}\prod_{\substack{A_1<p\le x/n_1\\p\mid m\text{~or~}p\in\mathscr{Q}}}\left(1-\frac{f(p)}{p}\right)\left(1-\frac{1}{p}\right)^{-1}\exp\left(\sum_{p\le x/n_1}\frac{f(p)-\nu(p,n_1;m)}{p}\right)\nonumber\\
&\ll\frac{y}{n_1\log2x}\prod_{\substack{p\le x/n_1\\p\nmid n_1,p\in\mathscr{Q}}}\left(1-\frac{1}{p}\right)^{-1}\exp\left(\sum_{p\le x/n_1}\frac{f(p)1_{p\nmid m,p\notin\mathscr{Q}}-\nu(p)1_{p\nmid n_1}}{p}\right)\nonumber\\
&\le\frac{y}{n_1\log2x}\prod_{p\in\mathscr{Q}}\left(1-\frac{1}{p}\right)^{-1}\exp\left(\sum_{p\le x/n_1}\frac{f(p)1_{p\nmid m,p\notin\mathscr{Q}}-\nu(p)1_{p\nmid n_1}}{p}\right),
\end{align}
where
\[\nu(p,n_1;m)\colonequals\begin{cases}
	~\nu(p),&\text{~~~~~if $p\nmid m$},\\
	~\nu(p)+1_{0\notin\mathscr{E}_p},&\text{~~~~~if $p\mid m$ and $p\nmid n_1$},\\
	~1,&\text{~~~~~if $p\mid n_1$}.
\end{cases}\]
By condition (i) and the constraint $\gcd(n_1,P_{\mathscr{E}})=1$, we see that the exponential in \eqref{eq:sumn_2} is 
\begin{align*}
&\ll \prod_{\substack{A_1<p\le x\\p\mid m\text{~or~}p\in\mathscr{Q}}}\left(1-\frac{f(p)}{p}\right)\prod_{\substack{p>k\\p\mid n_1}}\left(1-\frac{\nu(p)}{p}\right)^{-1}\exp\left(\sum_{p\le x}\frac{f(p)-\nu(p)}{p}\right)\\
&\le\prod_{\substack{p>A_1\\p\mid m,p\notin\mathscr{Q}}}\left(1-\frac{f(p)}{p}\right)\prod_{p\in\mathscr{Q}}\left(1-\frac{f(p)}{p}\right)\prod_{\substack{p>k\\p\mid m,p\notin\mathscr{Q}}}\left(1-\frac{\nu(p)}{p}\right)^{-1}\exp\left(\sum_{p\le x}\frac{f(p)-\nu(p)}{p}\right),
\end{align*}
which implies that
\begin{align*}
\sum_{\substack{n_1< z\\ \rad(n_1)\mid m\\ (n_1,P_{\mathscr{E}})=1}}f(mn_1)\sum_{\substack{n_1n_2\in\mathscr{N}\\ (n_2,m)=1}}f(n_2)&\ll C_{f,\mathscr{Q}}(m)\frac{y}{\log 2x}\exp\left(\sum_{p\le x}\frac{f(p)-\nu(p)}{p}\right)\sum_{\substack{n_1\ge1\\ \rad(n_1)\mid m}}\frac{f(mn_1)}{n_1}\\
&\le C_{f,\mathscr{Q}}(m)m\frac{y}{\log 2x}\exp\left(\sum_{p\le x}\frac{f(p)-\nu(p)}{p}\right)\prod_{p^v\parallel m}\sum_{\nu\ge v}\frac{f(p^{\nu})}{p^{\nu}}.
\end{align*}
On the other hand, if $z\le n_1\le x$, then we use the trivial upper bound
\[\sum_{\substack{ n_1n_2\in\mathscr{N}\\ (n_2,m)=1}}f(n_2)\ll \left(\frac{x}{z}\right)^{\epsilon}\sum_{(x-y)/n_1<n_2\le x/n_1}1\le \left(\frac{x}{z}\right)^{\epsilon}\frac{y}{z}=\left(\frac{x}{y}\right)^{2\epsilon/(1-\epsilon)}< x^{4\epsilon}\]
implied by condition (ii) and the fact that $0<\epsilon<1/2$. It follows that
\[\sum_{\substack{z\le n_1\le x\\ \rad(n_1)\mid m\\ (n_1,P_{\mathscr{E}})=1}}f(mn_1)\sum_{\substack{ n_1n_2\in\mathscr{N}\\ (n_2,m)=1}}f(n_2)\ll x^{4\epsilon}\sum_{\substack{n_1\le x\\ \rad(n_1)\mid m}}(mn_1)^{\epsilon}\le y^{\delta_0/2}\sum_{\substack{n_1\le x\\ \rad(n_1)\mid m}}1.\]
The sum in the last inequality above is maximized when $\rad(m)=p_1\cdots p_{\ell}$ is the product of the first $\ell$ primes $p_1<\cdots<p_{\ell}$ with $\ell=\omega(m)$. Since $\ell\le (\eta_0+o(1))\log x/\log\log x$, we have
\[\sum_{\substack{n_1\le x\\ \rad(n_1)\mid m}}1\le\Psi(x,p_{\ell})\le\frac{e^{O(\ell)}}{\ell!}\prod_{i=1}^{\ell}\frac{\log x}{\log p_{i}}\le e^{O(\ell)}\prod_{i=1}^{\ell}\frac{\log x}{\ell\log p_{i}}\]
by \cite[Theorem III.5.3]{Ten} and Stirling's formula, where $\Psi(x,p_{\ell})$ counts the number of $p_{\ell}$-smooth numbers $n\le x$. If $\ell\le(\alpha\delta_0/4)\log x/\log\log x$, then we have
\[\sum_{\substack{n_1\le x\\ \rad(n_1)\mid m}}1\ll e^{O(\ell)}(\log x)^{\ell}\le x^{\alpha\delta_0/2}<y^{\delta_0/2}\]
for sufficiently large $x$. If $\ell\ge(\alpha\delta_0/4)\log x/\log\log x$, then 
\[\prod_{i=1}^{\ell}\log p_i=\exp\left(\sum_{p\le p_{\ell}}\log\log p\right)=\exp\left((1+o(1))\frac{p_{\ell}\log\log p_{\ell}}{\log p_{\ell}}\right)=(\log\ell)^{(1+o(1))\ell},\]
which implies that
\[\sum_{\substack{n_1\le x\\ \rad(n_1)\mid m}}1\le e^{O(\ell)}\left(\frac{\log x}{\ell(\log \ell)^{1+o(1)}}\right)^{\ell}=e^{O(\ell)}(\log\log x)^{o(\ell)}= x^{o(1)}.\]
In either case, we can ensure that
\[\sum_{\substack{z\le n_1\le x\\ \rad(n_1)\mid m\\ (n_1,P_{\mathscr{E}})=1}}f(mn_1)\sum_{\substack{ n_1n_2\in\mathscr{N}\\ (n_2,m)=1}}f(n_2)\ll y^{\delta_0}\]
for sufficiently large $x$. This completes the proof of the lemma.
\end{proof}

\par We say that a collection $\{(a_i,b_i)\}_{1\le i\le k}$, where $a_i\in\Z\setminus\{0\}$ and $b_i\in\Z$ for every $1\le i\le k$, is {\it admissible} if $\rho(p)\colonequals \#\{n\in\Z/p\Z\colon (a_1n+b_1)\cdots(a_kn+b_k)\equiv0\psmod{p}\}<p$ for all primes $p$. The {\it singular series} $\mathfrak{S}$ associated to an admissible collection $\{(a_i,b_i)\}_{1\le i\le k}$ is defined by
\[\mathfrak{S}\colonequals\prod_{p}\left(1-\frac{\rho(p)}{p}\right)\left(1-\frac{1}{p}\right)^{-k}.\]
The following simple but useful corollary of Lemma \ref{lem:Pollf(mn)} is immediate.
\begin{lem}\label{lem:weightedHL}
Let $\alpha,\delta_0\in(0,1)$, $k\in\Z_{\ge0}$, $\eta_0,A_1>0$, $A_2\colon \R_{>0}\to\R_{>0}$, and $f\in\mathscr{M}_s(A_1,A_2)$. Then we have
\begin{align*}
&\hspace*{5mm}\sum_{x-y<n\le x}f(mn)1_{\PP}(a_1n+b_1)\cdots1_{\PP}(a_kn+b_k)\\
&\ll_{\alpha,k,\delta_0,\eta_0,A_1,A_2} C_{f,b}(m)m\frac{\mathfrak{S}y}{(\log 2x)^{k+1}}\exp\left(\sum_{p\le x}\frac{f(p)}{p}\right)\prod_{p^v\parallel m}\sum_{\nu\ge v}\frac{f(p^{\nu})}{p^{\nu}}+y^{\delta_0}\\
&\ll_kC_{f,b}(m)m\frac{y}{(\log 2x)^{k+1}}\exp\left(\sum_{p\le x}\frac{f(p)}{p}\right)\prod_{p\mid D}\left(1-\frac{1}{p}\right)^{\rho(p)-k}\prod_{p^v\parallel m}\sum_{\nu\ge v}\frac{f(p^{\nu})}{p^{\nu}}+y^{\delta_0}
\end{align*}
for all $x\ge 1$, $x^{\alpha}<y\le x$, $m\in\N\cap[1,x^{\eta_0}]$, and all admissible collections $\{(a_i,b_i)\}_{1\le i\le k}$ with 
\[D=D\left(\{(a_i,b_i)\}_{1\le i\le k}\right)\colonequals\prod_{i=1}^{k}a_i\prod_{1\le j<i}(a_ib_j-a_jb_i)\ne0,\]
where 
\[C_{f,b}(m)\colonequals\prod_{\substack{p>\max(A_1,k)\\p\mid m\\ p\nmid b_1\cdots b_k}}\left(1-\frac{f(p)}{p}\right)\left(1-\frac{\rho(p)}{p}\right)^{-1}\prod_{\substack{A_1<p\le x\\p\mid b_1\cdots b_k}}\left(1-\frac{f(p)}{p}\right)\left(1-\frac{1}{p}\right)^{-1}.\]
\end{lem}
\begin{proof}
The proof is essentially the same as that of the classical result for $f=1$ from sieve theory. In particular, the second asserted upper bound follows from the first and the well-known upper bound for $\mathfrak{S}$. To prove the first asserted upper bound, we assume that $x$ is sufficiently large depending on $\alpha,k,\delta_0,\eta_0,A_1,A_2$. In the proof, all the implicit constants will depend at most on these parameters. Let $P(y)\colonequals\prod_{p\le y^{\delta_0/2}}p$, $Q(n)\colonequals(a_1n+b_1)\cdots(a_kn+b_k)$, and
$\mathscr{N}\colonequals\{x-y<n\le x\colon\gcd(Q(n),P(y))=1\}$. Then 
\[\nu(p)=\begin{cases}
		~\rho(p)&\text{~~~~~if $p\le y^{\delta_0/2}$},\\
		~0&\text{~~~~~otherwise},
\end{cases}\]
and $\mathscr{Q}=\{A_1<p\le y^{\delta_0/2}\colon p\mid b_1\cdots b_k\}$ in this particular case. By condition (ii) we see that
\begin{align*}
\sum_{\substack{n\le x\\ 0<a_in+b_i\le y^{\delta_0/2}\text{~for some~}i}}f(mn)1_{\PP}(a_1n+b_1)\cdots1_{\PP}(a_kn+b_k)&\ll (mx)^{\frac{\alpha\delta_0}{2(1+\eta_0)}}\sum_{0<a_in+b_i\le y^{\delta_0/2}\text{~for some~}i}1\\
&\ll y^{\delta_0}.
\end{align*}
On the other hand, we have by Lemma \ref{lem:Pollf(mn)} and condition (i) that
\begin{align*}
\sum_{n\in\mathscr{N}}f(mn)&\ll C_{f,\mathscr{Q}}(m)m\frac{y}{\log 2x}\exp\left(\sum_{p\le x}\frac{f(p)-\nu(p)}{p}\right)\prod_{p^v\parallel m}\sum_{\nu\ge v}\frac{f(p^{\nu})}{p^{\nu}}+y^{\delta_0}\\
&\ll C_{f,b}(m)m\frac{y}{\log 2x}\exp\left(\sum_{p\le x}\frac{f(p)-\rho(p)}{p}\right)\prod_{p^v\parallel m}\sum_{\nu\ge v}\frac{f(p^{\nu})}{p^{\nu}}+y^{\delta_0}\\
&\ll C_{f,b}(m)m\frac{\mathfrak{S}y}{(\log 2x)^{k+1}}\exp\left(\sum_{p\le x}\frac{f(p)}{p}\right)\prod_{p^v\parallel m}\sum_{\nu\ge v}\frac{f(p^{\nu})}{p^{\nu}}+y^{\delta_0}.
\end{align*}
Adding up the two estimates above yields the first asserted upper bound.
\end{proof}

Lemma \ref{lem:weightedHL} allows us to obtain instantly a generalization of Proposition \ref{prop:philcm} for sub-multiplicative functions satisfying the conditions (i) and (ii). Before proving such a generalization, we establish the following technical lemma.
\begin{lem}\label{lem:fphi(av-bu)}
Let $a,b\in\Z\setminus\{0\}$, $A_1>0$, $A_2\colon \R_{>0}\to\R_{>0}$, and $f\in\mathscr{M}_s(A_1,A_2)$. Then 
\[\sum_{\substack{uv\le x\\ (u,av)=(v,bu)=1\\ av\ne bu}}f(uv)\cdot\frac{|av-bu|}{\varphi(|av-bu|)}\ll_{a,b,A_1,A_2}\frac{x}{\log x}\exp\left(2\sum_{p\le x}\frac{f(p)}{p}\right)\]
for all $x\ge2$.
\end{lem}
\begin{proof}
It is sufficient to consider the case where $x$ is sufficiently large depending on $a,b,A_1,A_2$. All the implicit constants appearing in the proof will depend at most on these parameters. In view of the fact that $\id/\varphi$ is completely sub-multiplicative, we may assume that $\gcd(a,b)=1$. Using the inequalities $n/\varphi(n)\ll\sigma(n)/n$ and 
\[\frac{\sigma(n)}{n}\le2\sum_{\substack{d\mid n\\ d\le\sqrt{n}}}\frac{1}{d},\]
we find that
\[\sum_{\substack{uv\le x\\ (u,av)=(v,bu)=1\\ av\ne bu}}f(uv)\cdot\frac{|av-bu|}{\varphi(|av-bu|)}\ll\sum_{d\ge1}\frac{1}{d}\sum_{\substack{uv\le x\\ av\ne bu,(av,bu)=1\\ av\equiv bu\psmod{d}\\\sqrt{|av-bu|}\ge d}}f(uv).\]
It suffices to consider the case $u\le v$, since the antithetical case $v\le u$ is similar. Now the conditions $u\le v$ and $d\le\sqrt{|av-bu|}$ imply that $d\ll\sqrt{x}$ and $v\gg d^2$. Since $\gcd(av,bu)=1$, we know that $\gcd(abu,d)=1$ and the congruence $av\equiv bu\psmod{d}$ is thus equivalent to $v\equiv \bar{a}_dbu\psmod{d}$, where $\bar{a}_d$ is the inverse of $a$ in $(\Z/d\Z)^{\times}$. By Shiu's theorem \cite[Theorem 1]{Shiu}, we obtain
\[\sum_{\substack{u\le v\le x/u\\ v\gg d^2,av\ne bu\\ (v,bu)=1\\ av\equiv bu\psmod{d}}}f(v)\le\sum_{\substack{v\le x/u\\ v\equiv \bar{a}_dbu\psmod{d}}}f(v)\ll\frac{1}{u\varphi(d)}\cdot\frac{x}{\log x}\exp\left(\sum_{p\le x, p\nmid d}\frac{f(p)}{p}\right)\]
whenever $x/u\gg d^2$, which actually follows from $d^2\ll v\le x/u$. Hence, we have
\begin{align*}
\sum_{\substack{uv\le x\\ (u,av)=(v,bu)=1\\ av\ne bu}}f(uv)\cdot\frac{|av-bu|}{\varphi(|av-bu|)}&\ll\frac{x}{\log x}\exp\left(\sum_{p\le x}\frac{f(p)}{p}\right)\sum_{d\ge1}\frac{1}{d\varphi(d)}\sum_{u\le\sqrt{x}}\frac{f(u)}{u}\\
&\ll\frac{x}{\log x}\exp\left(2\sum_{p\le x}\frac{f(p)}{p}\right),
\end{align*}
completing the proof of the lemma.
\end{proof}

We are now in a position to prove the following generalization of Proposition \ref{prop:philcm}.
\begin{thm}\label{thm:f([p-a,q-a])}
Let $a,b\in\Z\setminus\{0\}$, $A_1>0$, $A_2\colon \R_{>0}\to\R_{>0}$, and $f\in\mathscr{M}_s(A_1,A_2)$. Then 
\[\sum_{\substack{[p-a,q-b]\le x\\p>a,q>b}}f([p-a,q-b])\ll_{a,b,A_1,A_2}E_f(x)\frac{x}{(\log x)^2}\int_{1}^{x}\frac{E_f(t)^2}{t(\log 2t)^2}\,dt\]
for all $x\ge 2$, where 
\[E_f(x)\colonequals\exp\left(\sum_{p\le x}\frac{f(p)}{p}\right).\]
\end{thm}
\begin{proof}
The proof is analogous to that of Proposition \ref{prop:philcm}. Before embarking on the proof, we introduce the multiplicative function
 \[\tilde{f}(n)\colonequals n\prod_{p^v\parallel n}\sum_{\nu\ge v}\frac{f(p^{\nu})}{p^{\nu}},\]
and observe that $\tilde{f}(p^v)=f(p^v)+O(1/p)$ for all prime powers $p^{\nu}$ and that $\tilde{f}\in\mathscr{M}_s\big(\widetilde{A}_1,\widetilde{A}_2\big)$ for some $\widetilde{A}_1>0$ and $\widetilde{A}_2\colon\R_{>0}\to\R_{>0}$ depending on $A_1$ and $A_2$. 
\par In what follows, we assume that $x$ is sufficiently large depending on $a,b,A_1,A_2$, and all the implicit constants appearing below depend at most on these parameters. Let $p-a=uw$ and $q-b=vw$ with $w=\gcd(p-a,q-b)$, so that $\gcd(u,v)=1$ and $[p-a,q-b]=uvw$. Then the sum to be estimated becomes
\begin{equation}\label{eq:fuvw1}
\sum_{\substack{uvw\le x\\ (u,v)=1}}f(uvw)1_{\PP}(uw+a)1_{\PP}(vw+b).
\end{equation}
\par As before, we handle the boundary cases first. Let us start with the contribution to \eqref{eq:fuvw1} from the terms with $\gcd(w,a)>1$ or $\gcd(w,b)>1$. By symmetry, it suffices to consider the contribution from the terms with $\gcd(w,a)=d>1$. Clearly, these terms vanish unless $uw+a=d$ is prime. Assuming that $d-a\ge1$, we see that the contribution from these terms is at most
\[\sum_{w\mid d-a}\sum_{v\le x/(d-a)}f((d-a)v)1_{\PP}(vw+b)\ll E_f(x)\frac{x}{(\log x)^2}\]
by Lemma \ref{lem:weightedHL} and condition (i). Similarly, one shows that this is also an upper bound for the contribution from the terms with $\gcd(u,a)>1$ or $\gcd(v,b)>1$. Hence, we may assume now that $\gcd(uw,a)=\gcd(vw,b)=1$. The remaining boundary case is $av=bu$. However, since $\gcd(u,a)=\gcd(v,b)=\gcd(u,v)=1$, this is impossible unless $a=b$ and $u=v=1$. Thus, the contribution from the terms in this case is also at most
\[\sum_{w\le x}f(w)1_{\PP}(w+a)\ll E_f(x)\frac{x}{(\log x)^2}.\]
This settles all the boundary cases.
\par From now on, we shall concentrate on the case where $\gcd(uw,a)=\gcd(vw,b)=1$ and $av\ne bu$. Fix any $\epsilon_0\in(0,1/3)$, whose exact value will be of little importance. We start by estimating the contribution to \eqref{eq:fuvw1} from the terms in this case with $uv\le x^{2/3}$. By Lemma \ref{lem:weightedHL} with $\delta_0=3\epsilon_0$ we see that 
\begin{align*}
&\hspace*{5mm}\sum_{\substack{w\le x/uv\\ (w,ab)=1}}f(uvw)1_{\PP}(uw+a)1_{\PP}(vw+b)\\
&\ll\frac{\tilde{f}(uv)}{uv}\left(\frac{uv}{\varphi(uv)}\right)^2\frac{|av-bu|}{\varphi(|av-bu|)}\cdot E_f(x)\frac{x}{(\log x)^3}+\left(\frac{x}{uv}\right)^{3\epsilon_0}.
\end{align*}
Summing this over $uv\le x^{2/3}$ with $\gcd(u,av)=\gcd(v,bu)=1$ yields the contribution
\[E_f(x)\frac{x}{(\log x)^3}\sum_{\substack{uv\le x^{2/3}\\ (u,av)=(v,bu)=1\\ av\ne bu}}\frac{\tilde{f}(uv)}{uv}\left(\frac{uv}{\varphi(uv)}\right)^2\frac{|av-bu|}{\varphi(|av-bu|)}+x^{\frac{2}{3}+\epsilon_0}\log x,\]
which, by Lemma \ref{lem:fphi(av-bu)} with $f$ replaced by $\tilde{f}(\id/\varphi)^2$ and partial summation, is 
\[\ll E_f(x)\frac{x}{(\log x)^3}\left(\frac{E_f(x)^2}{\log x}+\int_{1}^{x}\frac{E_f(t)^2}{t\log 2t}\,dt\right).\]
\par It remains to consider the contribution to \eqref{eq:fuvw1} from the terms with $uv>x^{2/3}$. Without loss of generality, we may assume $u\le v$. In this case, we have $w<x^{1/3}$ and $u\le\sqrt{x/w}$, so that $uw\le\sqrt{wx}<x^{2/3}$. By Lemma \ref{lem:weightedHL} with $\delta_0=3\epsilon_0$ we have
\[\sum_{v\le x/uw}f(uvw)1_{\PP}(vw+b)\ll\frac{\tilde{f}(uw)}{uw}\cdot\frac{uw}{\varphi(uw)}\cdot\frac{w}{\varphi(w)}\cdot E_f(x)\frac{x}{(\log x)^2}+\left(\frac{x}{uw}\right)^{3\epsilon_0}.\]
It follows that the contribution from the present case is
\begin{align*}
&\ll E_f(x)\frac{x}{(\log x)^2}\sum_{uw< x^{2/3}}\frac{\tilde{f}(uw)}{uw}\left(\frac{uw}{\varphi(uw)}\right)^21_{\PP}(uw+a)+x^{\frac{2}{3}+\epsilon_0}\\
&\ll E_f(x)\frac{x}{(\log x)^2}\sum_{n< x}\frac{\tilde{f}(n)\tau(n)}{n}\left(\frac{n}{\varphi(n)}\right)^21_{\PP}(n+a)+x^{\frac{2}{3}+\epsilon_0}.
\end{align*}
By Lemma \ref{lem:weightedHL}, we have
\[\sum_{n< x}\tilde{f}(n)\tau(n)\left(\frac{n}{\varphi(n)}\right)^21_{\PP}(n+a)\ll E_f(x)^2\frac{x}{(\log x)^2}.\]
By partial summation, we see that the contribution from the present case is
\[\ll E_f(x)\frac{x}{(\log x)^2}\left(\frac{E_f(x)^2}{(\log x)^2}+\int_{1}^{x}\frac{E_f(t)^2}{t(\log 2t)^2}\,dt\right).\]
Since 
\begin{align*}
\int_{1}^{x}\frac{E_f(t)^2}{t(\log 2t)^2}\,dt&\ge\int_{1}^{x}\frac{1}{t(\log 2t)^2}\,dt\gg1,\\
\int_{1}^{x}\frac{E_f(t)^2}{t(\log 2t)^2}\,dt&\gg\frac{1}{\log x}\int_{1}^{x}\frac{E_f(t)^2}{t\log 2t}\,dt,\\
\int_{1}^{x}\frac{E_f(t)^2}{t(\log 2t)^2}\,dt&\gg E_f(x)^2\int_{\sqrt{x}}^{x}\frac{1}{t(\log 2t)^2}\,dt\gg\frac{E_f(x)^2}{\log x},
\end{align*}
the desired upper bound follows upon collecting the estimates for the contributions to \eqref{eq:fuvw1} from all of the cases above.
\end{proof}

The following special case of Theorem \ref{thm:f([p-a,q-a])} when $f=\tau_{\kappa}$ with $\kappa\ge0$ is immediate.
\begin{cor}\label{cor:tau_k([p-a,q-b])}
Let $a,b\in\Z\setminus\{0\}$ and $\kappa\ge0$. Then 
\[\sum_{\substack{[p-a,q-b]\le x\\p>a,q>b}}\tau_{\kappa}([p-a,q-b])\ll_{a,b,\kappa}x(\log x)^{\max(\kappa-2,3\kappa-3)}(\log\log x)^{1_{\kappa=1/2}}\]
for all $x\ge3$.
\end{cor}

Corollary \ref{cor:tau_k([p-a,q-b])} shows that if we define the natural probability measure induced by $\tau_{\kappa}$ on $\N\cap[1,x]$ by
\[\text{Prob}(n=n_0)=\frac{\tau_{\kappa}(n_0)}{\sum_{m\le x}\tau_{\kappa}(m)}\] 
for every $n_0\in\N\cap[1,x]$, then with respect to this probability measure, the expected number of representations of $n\in\N\cap[1,x]$ of the form $n=[p-a,q-b]$ with some primes $p>a$ and $q>b$ is $\ll_{a,b,\kappa}(\log x)^{\max(-1,2\kappa-2)}(\log\log x)^{1_{\kappa=1/2}}$.
\par More generally, it may be of interest to estimate
\[\sum_{\substack{[p_1-a_1,...,p_k-a_k]\le x\\p_i>a_i}}f([p_1-a_1,...,p_k-a_k]),\]
where $k\ge2$, $\boldsymbol{a}=(a_1,...,a_k)\in(\Z\setminus\{0\})^{k}$, and $f\in\mathscr{M}_s(A_1,A_2)$. Perhaps we have
\[\sum_{\substack{[p_1-a_1,...,p_k-a_k]\le x\\p_i>a_i}}f([p_1-a_1,...,p_k-a_k])\ll_{\boldsymbol{a},A_1,A_2}E_f(x)\frac{x}{(\log x)^2}\int_{1}^{x}\frac{E_f(t)^{2^k-2}}{t(\log 2t)^k}\,dt\]
for all $x\ge2$. In particular, if $a_1=\cdots=a_k=a$, then the above estimate applied to $f=1$ would imply 
\[\frac{1}{x}\sum_{n\le x}\omega^*_a(n)^k\ll_{k,a} (\log x)^{2^k-k-1},\]
yielding an upper bound of the conjectured order of magnitude for $M_k(x)$.
\medskip
\section{A conjecture on the second moment of $\omega^*$}\label{S:2Momega*}
In this section, we discuss a possible asymptotic formula for the second moment of $\omega^*$ of the form
\[M_2(x)=\frac{1}{x}\sum_{n\le x}\omega^*(n)^2\sim C\log x,\]
where $C>0$ is some constant. We define the related quantity 
\[S_2(x)\colonequals\frac{1}{x}\cdot\#\{(p,q)\colon [p-1,q-1]\le x\}=\frac{1}{x}\sum_{n\le x}\beta_1(n),\]
where $\beta_1(n)$ was first introduced in the proof of Theorem \ref{thm:2mshiftomega_a^*}. It can be shown by partial summation that $S_2(x)\sim C$ implies $M_2(x)\sim C\log x$. Based on the Elliott--Halberstam conjecture, Ding, Guo, and Zhang \cite{DGZ} argued that $C=2\zeta(2)\zeta(3)/\zeta(6)=315\zeta(3)/\pi^4\approx3.88719$. However, Pomerance and the author \cite{FP} discovered an error in \cite[Equation (4.8)]{MM} on which their heuristic argument is based. As discussed in  \cite[Section 2]{FP}, numerical computations seem to suggest that $C\approx 3.2$. Here we propose a crude heuristic argument leading us to the following conjecture.
\begin{conj}\label{conj:S_2}
We have
\begin{align*}
S_2(x)&\sim\frac{105}{4\pi^2}\zeta(3),\\
M_2(x)&\sim\frac{105}{4\pi^2}\zeta(3)\log x.
\end{align*}
\end{conj}
In other words, we conjecture that
\[C=\frac{105}{4\pi^2}\zeta(3)=3.1970879911....\]
In contrast to the constant $2\zeta(2)\zeta(3)/\zeta(6)$ obtained in \cite{DGZ}, the constant that we derived is actually equal to $\zeta(2)^2\zeta(3)/\zeta(6)$. Our heuristic argument relies on the following variation of the uniform Bateman--Horn conjecture for linear polynomials: 
\begin{equation}\label{eq:HLab}
\sum_{n\le x}1_{\PP}(an+1)1_{\PP}(bn+1)=\mathfrak{S}(a,b)\frac{x}{\log(2ax)\log(2bx)}+E(x;a,b)
\end{equation}
for all $a,b\in\N$ with $ab\le x$, $a\ne b$ and $\gcd(a,b)=1$, where the error $E(x;a,b)$ is expected to be negligible on average compared to the main term. Here
\[\mathfrak{S}(a,b)\colonequals \prod_{p}\left(1-\frac{\nu_p(a,b)}{p}\right)\left(1-\frac{1}{p}\right)^{-2}=2C_2\prod_{\substack{p>2\\p\mid ab(a-b)}}\left(1+\frac{1}{p-2}\right)\]
is the singular series corresponding to the linear polynomials $an+1$ and $bn+1$, where $\nu_p(a,b)\colonequals\#\{n\in\Z/p\Z\colon (an+1)(bn+1)\equiv0\psmod{p}\}$, and 
\[C_2\colonequals\prod_{p>2}\left(1-\frac{1}{(p-1)^2}\right)\]
is the twin prime constant. Using \eqref{eq:HLab} we expect to have
\begin{align*}
S_2(x)&=\frac{\pi(x+1)}{x}+\frac{1}{x}\sum_{\substack{ab\le x\\a\ne b\\(a,b)=1}}\sum_{n\le x/ab}1_{\PP}(an+1)1_{\PP}(bn+1)\\
&\sim\sum_{\substack{ab\le x\\a\ne b\\(a,b)=1}}\frac{\mathfrak{S}(a,b)}{ab\log(2x/a)\log(2x/b)}\\
&=\sum_{k,l\ge0}\frac{1}{(\log2x)^{k+l+2}}\sum_{\substack{ab\le x\\a\ne b\\(a,b)=1}}\frac{\mathfrak{S}(a,b)}{ab}(\log a)^k(\log b)^l.
\end{align*}
Instead of dealing with the infinite product in the original definition of $\mathfrak{S}(a,b)$ given above, we rewrite it as an infinite series. For every integer $\nu\ge0$ and every prime $p$, we write
\begin{equation}\label{eq:s(p,nu)}
s(p,\nu)\colonequals\left(1-\frac{\nu}{p}\right)\left(1-\frac{1}{p}\right)^{-2}-1=\frac{(2-\nu)p-1}{(p-1)^2},
\end{equation}
and put
\[s_{a,b}(n)\colonequals \prod_{p\mid n}s(p,\nu_p(a,b))\]
if $n$ is square-free and $s_n(a,b)\colonequals 0$ otherwise. Then we have
\[\mathfrak{S}(a,b)=\frac{ab}{\varphi(ab)}\sum_{(n,ab)=1}s_{a,b}(n)=\frac{ab}{\varphi(ab)}\left(\sum_{\substack{n\le y\\(n,ab)=1}}s_{a,b}(n)+O\left(\frac{3^{\omega(|a-b|)}}{y}\right)\right)\]
for $y\ge2$ and $a,b\in\N$ with $a\ne b$ and $\gcd(a,b)=1$, where the second equality follows from
\begin{align*}
\sum_{\substack{n>y\\(n,ab)=1}}|s_{a,b}(n)|&\le\sum_{\substack{n>y\\(n,ab)=1}}\mu(n)^2\frac{\varphi((n,a-b))}{\varphi(n)^2}\\
&=\sum_{d\mid a-b}\frac{\mu(d)^2}{\varphi(d)}\sum_{\substack{m>y/d\\(m,ab(a-b))=1}}\frac{\mu(m)^2}{\varphi(m)^2}\\
&\ll\frac{1}{y}\sum_{d\mid a-b}\frac{\mu(d)^2d}{\varphi(d)}\le\frac{3^{\omega(|a-b|)}}{y}
\end{align*} 
by \eqref{eq:s(p,nu)}. Thus, we expect
\[S_2(x)\sim\sum_{k,l\ge0}\frac{1}{(\log2x)^{k+l+2}}\sum_{\substack{ab\le x\\a\ne b\\(a,b)=1}}\sum_{\substack{n\le y\\(n,ab)=1}}\frac{s_{a,b}(n)}{\varphi(a)\varphi(b)}(\log a)^k(\log b)^l\]
for some suitable choice of $y$ in terms of $x$. Given any square-free $n\in\N$ and any sequence $\{\nu_p\}_{p\mid n}\subseteq\{1,2\}^{\omega(n)}$, let
\[\mathcal{A}_n\colonequals\{(a_n,b_n)\in\left(\Z/n\Z\right)^{\times}\times \left(\Z/n\Z\right)^{\times}\colon \nu_p(a_n,b_n)=\nu_p\text{~for all~}p\mid n\}.\]
Then the inner sum over $ab\le x$ above becomes
\[\sum_{n\le y}\mu(n)^2\sum_{\{\nu_p\}_{p\mid n}\subseteq\{1,2\}^{\omega(n)}}\prod_{p\mid n}s(p,\nu_p(a,b))\sum_{(a_n,b_n)\in\mathcal{A}_n}\sum_{\substack{a\le x\\ a\equiv a_n\psmod{n}}}\frac{(\log a)^k}{\varphi(a)}\sum_{\substack{b\le x/a\\ b\ne a\\(b,a)=1\\ b\equiv b_n\psmod{n}}}\frac{(\log b)^l}{\varphi(b)}.\]
Since
\begin{align*}
\sum_{\substack{b\le z\\ b\ne a\\(b,a)=1\\ b\equiv b_n\psmod{n}}}\frac{b}{\varphi(b)}&=\sum_{\substack{d\le z\\(d,an)=1}}\frac{\mu(d)^2}{\varphi(d)}\sum_{\substack{b\le z/d\\ bd\ne a\\(b,a)=1\\ bd\equiv b_n\psmod{n}}}1\\
&=\varphi(a)\sum_{\substack{d\le z\\(d,an)=1}}\frac{\mu(d)^2}{\varphi(d)}\left(\frac{z}{dan}+O(1)\right)\\
&=B_0\frac{\varphi(a)}{an}z\prod_{p\mid an}\left(1+\frac{1}{p(p-1)}\right)^{-1}+O(\varphi(a)\log z),
\end{align*}
where 
\[B_0\colonequals\prod_{p}\left(1+\frac{1}{p(p-1)}\right)=\frac{\zeta(2)\zeta(3)}{\zeta(6)},\]
we expect
\[\sum_{\substack{b\le x/a\\ b\ne a\\(b,a)=1\\ b\equiv b_n\psmod{n}}}\frac{(\log b)^l}{\varphi(b)}\sim\frac{B_0}{l+1}\cdot\frac{\varphi(a)g(a)g(n)}{an}\left(\log\frac{x}{a}\right)^{l+1},\]
where 
\[g(n)\colonequals\prod_{p\mid n}\left(1+\frac{1}{p(p-1)}\right)^{-1}.\]
Thus, we expect
\[\sum_{\substack{a\le x\\ a\equiv a_n\psmod{n}}}\frac{(\log a)^k}{\varphi(a)}\sum_{\substack{b\le x/a\\ b\ne a\\(b,a)=1\\ b\equiv b_n\psmod{n}}}\frac{(\log b)^l}{\varphi(b)}\sim\frac{B_0}{l+1}\cdot\frac{g(n)}{n}\sum_{\substack{a\le x\\ a\equiv a_n\psmod{n}}}\frac{g(a)}{a}(\log a)^k\left(\log\frac{x}{a}\right)^{l+1}.\]
Since $g=1\ast g_1$ with
\[g_1(n)=\mu(n)\prod_{p\mid n}\frac{1}{p^2-p+1},\]
we have
\[\sum_{\substack{a\le x\\ a\equiv a_n\psmod{n}}}g(a)=\sum_{\substack{d\le x\\(d,n)=1}}g_1(d)\left(\frac{x}{dn}+O(1)\right)=A_0x\prod_{p\mid n}\frac{p^2-p+1}{(p-1)(p^2+1)}+O(1),\]
where 
\[A_0\colonequals\prod_{p}\left(1-\frac{1}{p(p^2-p+1)}\right)=\frac{\zeta(6)}{\zeta(3)\zeta(4)}.\]
Hence, we expect
\[\sum_{\substack{a\le x\\ a\equiv a_n\psmod{n}}}\frac{g(a)}{a}(\log a)^k\sim\frac{A_0}{k+1}(\log x)^{k+1}\prod_{p\mid n}\frac{p^2-p+1}{(p-1)(p^2+1)}\]
and 
\begin{align*}
\sum_{\substack{a\le x\\ a\equiv a_n\psmod{n}}}\frac{g(a)}{a}(\log a)^k\left(\log\frac{x}{a}\right)^{l+1}&\sim\frac{A_0(l+1)}{k+1}\prod_{p\mid n}\frac{p^2-p+1}{(p-1)(p^2+1)}\int_{1}^{x}(\log t)^{k+1}\left(\log\frac{x}{t}\right)^{l}\frac{dt}{t}\\
&=\frac{A_0(l+1)B(k+1,l+1)}{k+l+2}(\log x)^{k+l+2}\prod_{p\mid n}\frac{p^2-p+1}{(p-1)(p^2+1)},
\end{align*}
where we have made the substitution $t=x^u$ and taken advantage of the beta function
\[B(\alpha,\beta)\colonequals\int_{0}^{1}u^{\alpha-1}(1-u)^{\beta-1}\,du\]
for any $\alpha,\beta>0$. Therefore, we expect
\[\sum_{\substack{a\le x\\ a\equiv a_n\psmod{n}}}\frac{(\log a)^k}{\varphi(a)}\sum_{\substack{b\le x/a\\ b\ne a\\(b,a)=1\\ b\equiv b_n\psmod{n}}}\frac{(\log b)^l}{\varphi(b)}\sim\frac{A_0B_0B(k+1,l+1)}{k+l+2}\eta(n)(\log x)^{k+l+2},\]
where
\[\eta(n)\colonequals \frac{g(n)}{n}\prod_{p\mid n}\frac{p^2-p+1}{(p-1)(p^2+1)}=\prod_{p\mid n}\frac{1}{p^2+1}.\]
At this point, we anticipate
\[S_2(x)\sim A_0B_0\sum_{k,l\ge1}\frac{B(k,l)}{k+l}\sum_{n\le y}\mu(n)^2\eta(n)\sum_{\{\nu_p\}_{p\mid n}\subseteq\{1,2\}^{\omega(n)}}\prod_{p\mid n}s(p,\nu_p(a,b))\sum_{(a_n,b_n)\in\mathcal{A}_n}1.\]
Since 
\[\sum_{(a_n,b_n)\in\mathcal{A}_n}1=\prod_{\substack{p\mid n\\\nu_p=1}}(p-1)\prod_{\substack{p\mid n\\\nu_p=2}}(p-1)(p-2),\]
it follows from \eqref{eq:s(p,nu)} that
\[\prod_{p\mid n}s(p,\nu_p(a,b))\sum_{(a_n,b_n)\in\mathcal{A}_n}1=\prod_{\substack{p\mid n\\\nu_p=2}}\left(\frac{1}{p-1}-1\right),\]
whence
\[\sum_{\{\nu_p\}_{p\mid n}\subseteq\{1,2\}^{\omega(n)}}\prod_{p\mid n}s(p,\nu_p(a,b))\sum_{(a_n,b_n)\in\mathcal{A}_n}1=\prod_{p\mid n}\left(1+\left(\frac{1}{p-1}-1\right)\right)=\frac{1}{\varphi(n)}\]
for square-free $n\in\N$. On the other hand, it is easy to see that
\begin{align*}
\sum_{k,l\ge1}\frac{B(k,l)}{k+l}&=\int_{0}^{1}\int_{0}^{1}\sum_{k,l\ge1}t^{k+l-1}u^{k-1}(1-u)^{l-1}\,dtdu\\
&=\int_{0}^{1}\int_{0}^{1}\frac{t}{1-ut}\cdot\frac{1}{1-(1-u)t}\,dudt\\
&=\int_{0}^{1}\int_{0}^{t}\frac{1}{1-z}\cdot\frac{1}{1-t+z}\,dzdt\\
&=\int_{0}^{1}\int_{t}^{1}\frac{1}{1-z}\cdot\frac{1}{1-t+z}\,dtdz\\
&=-\int_{0}^{1}\frac{\log z}{1-z}\,dz\\
&=\sum_{n\ge1}\frac{1}{n}\int_{0}^{1}(1-z)^{n-1}\,dz=\zeta(2).
\end{align*}
Hence, we anticipate
\[S_2(x)\sim A_0B_0\zeta(2)\sum_{n\ge1}\mu(n)^2\frac{\eta(n)}{\varphi(n)}=\frac{\zeta(2)^2}{\zeta(4)}\sum_{n\ge1}\mu(n)^2\frac{\eta(n)}{\varphi(n)}.\]
But 
\[\sum_{n\ge1}\mu(n)^2\frac{\eta(n)}{\varphi(n)}=\prod_{p}\left(1+\frac{1}{(p-1)(p^2+1)}\right)=\prod_{p}\frac{1-p^{-1}+p^{-2}}{(1-p^{-1})(1+p^{-2})}=\frac{\zeta(3)\zeta(4)}{\zeta(6)}.\]
Therefore, we expect to have
\[S_2(x)\sim \frac{\zeta(2)^2\zeta(3)}{\zeta(6)}=\frac{105}{4\pi^2}\zeta(3),\]
which is Conjecture \ref{conj:S_2}. It is possible to obtain heuristics for the $k$th moment $M_k(x)$ for all $k\ge2$ by generalizing the above argument.
\par Admittedly, we have ignored in the above argument all the possible error terms in the asymptotics, especially when handling the technicalities arising from various sums over arithmetic progressions where the moduli in consideration may be very large. Nonetheless, it is probable that the total contributions from these error terms do not affect the above asymptotic for $S_2(x)$. Furthermore, numerical computations seem to support our conjecture. 
\par Despite the incorrect value of the constant $C$ obtained in \cite{DGZ}, the argument given there does seem to suggest that
\[\sum_{p,q\le x}\frac{1}{[p-1,q-1]}\sim\frac{2\zeta(2)\zeta(3)}{\zeta(6)}\log x.\]
It is observed in \cite[Section 8]{FP} that 
\[\frac{1}{\log x}\sum_{\substack{p,q\le x\\ [p-1,q-1]>x}}\frac{1}{[p-1,q-1]}\approx 0.69,\]
based on numerical computations. Note that 
\[\frac{2\zeta(2)\zeta(3)}{\zeta(6)}-\frac{\zeta(2)^2\zeta(3)}{\zeta(6)}=0.6901048825....\]
This provides further evidence for the truthfulness of Conjecture \ref{conj:S_2}.
\par Finally, it has also been conjectured in a footnote in \cite{FP} that 
\begin{equation}\label{eq:difM_2/S_2}
S_2(x)^{-1}\left(\sum_{[p-1,q-1]\le x}\frac{1}{[p-1,q-1]}-M_2(x)\right)\sim 1-\gamma,
\end{equation}
where $\gamma=0.57722...$ is the Euler--Mascheroni constant. Accepting the conjecture that $S_2(x)\sim C$ for some constant $C>0$, one would not be surprised if this conjecture turns out to be true as well. Indeed, we have
\begin{equation}\label{eq:difM_2}
\sum_{[p-1,q-1]\le x}\frac{1}{[p-1,q-1]}-M_2(x)=\frac{1}{x}\sum_{[p-1,q-1]\le x}\left\{\frac{x}{[p-1,q-1]}\right\}=\frac{1}{x}\sum_{n\le x}\beta_1(n)\left\{\frac{x}{n}\right\}.
\end{equation}
By partial summation, we expect that the last expression in \eqref{eq:difM_2} is equal to
\[\frac{1}{x}\int_{1^-}^{x}\left\{\frac{x}{t}\right\}\,d(tS_2(t))\sim \frac{C}{x}\int_{1^-}^{x}\left\{\frac{x}{t}\right\}\,dt\sim C\int_{1}^{\infty}\frac{\{t\}}{t^2}\,dt=C(1-\gamma).\]
As intuitive as this argument is, making the first asymptotic above rigorous would require a careful treatment of the technicalities resulting from the Stieltjes integral involving the error term in $tS_2(t)$. Rather than try to overcome these technicalities, we give an argument which is based on a general principle and circumvents such technicalities.
\begin{prop}\label{prop:difM_2}
The conjectured asymptotic formula \eqref{eq:difM_2/S_2} holds if there exists an absolute constant $C>0$ such that $S_2(x)\sim C$ as $x\to\infty$.
\end{prop}
\begin{proof}
The proof is simple. Let $a_n\colonequals \beta_1(n)-C$ for $n\in\N$. Then 
\begin{align*}
\sum_{n\le x}a_n&=o(x),\\
\sum_{n\le x}|a_n|&=O(x).
\end{align*}
Applying Axer's theorem \cite[Theorem 8.1]{MV} to $\{a_n\}_{n\ge1}$ and $F(x)=\{x\}$, we have
\[\sum_{n\le x}a_n\left\{\frac{x}{n}\right\}=o(x).\]
Since
\[\sum_{n\le x}\left\{\frac{x}{n}\right\}=x\sum_{n\le x}\frac{1}{n}-\sum_{n\le x}\left\lfloor\frac{x}{n}\right\rfloor=x\sum_{n\le x}\frac{1}{n}-\sum_{n\le x}\tau(n)=(1-\gamma)x+O\left(\sqrt{x}\right),\]
we obtain
\[\frac{1}{x}\sum_{n\le x}\beta_1(n)\left\{\frac{x}{n}\right\}=\frac{1}{x}\sum_{n\le x}a_n\left\{\frac{x}{n}\right\}+\frac{C}{x}\sum_{n\le x}\left\{\frac{x}{n}\right\}=C(1-\gamma)+o(1),\]
which, in view of \eqref{eq:difM_2}, completes the proof of the proposition.
\end{proof}
\medskip
\section{Concluding remarks}\label{S:CR}
Given any arithmetic function $f\colon\N\to\N$ and any nonempty subset $S\subseteq\N$, we define the {\it $k$-level set} $\mathcal{L}_k(f,S)$ of the restriction $f|_{S}$ by
\[\mathcal{L}_k(f,S)\colonequals\{n\in S\colon f(n)=k\}\]
for each $k\in\N$. The natural density $\delta_k(f,S)$ of $\mathcal{L}_k(f,S)$ relative to $S$ is then given by
\[\delta_k(f,S)\colonequals\lim_{x\to\infty}\frac{\#(\mathcal{L}_k(f,S)\cap[1,x])}{\#(S\cap[1,x])},\]
provided that this limit exists. Despite the similarities between $\omega^*$ and $\tau$ suggested by their maximal orders and moments, the natures of the level sets $\mathcal{L}_k(\omega^*,\N)$ and $\mathcal{L}_k(\tau,\N)$ are quite different. It is well-known that both $\mathcal{L}_k(\tau,\N)$ and $\mathcal{L}_k(\omega,\N)$ have natural density 0 for every $k\in\N$. However, Pomerance and the author \cite[Theorem 2]{FP} recently proved that $\delta_k(\omega^*,\N)>0$ for every $k\in\N$ and that $\sum_{k\ge1}\delta_k(\omega^*,\N)=1$.
\par It is also of interest to examine these densities with $\N$ replaced by the set $\PP_b\colonequals\{p-b\colon p\in\PP\cap(b,\infty)\}$, where $b\in\Z\setminus\{0\}$. For the sake of simplicity, we stick with $\omega^*$ rather than pursue the more sophisticated function $\omega_a^*$ studied in Section \ref{S:omega_a*(p-b)}. A general result of the Erd\H{o}s--Kac type due to Halberstam \cite[Theorem 3]{Hal1} implies that $\delta_k(\omega,\PP_b)=\delta_k(\tau,\PP_b)=0$ for every $k\in\N$. The situation on $\delta_k(\omega^*,\PP_b)$ is somewhat complicated. If $b$ is even, then we clearly have $\delta_1(\omega^*,\PP_b)=1$ and $\delta_k(\omega^*,\PP_b)=0$ for all $k\ge2$. If $b$ is odd, then $\delta_1(\omega^*,\PP_b)=0$. In the special case $b=1$, we also have $\delta_2(\omega^*,\PP_b)=0$. Nevertheless, we expect that the densities $\delta_k(\omega^*,\PP_b)$ all exist and add up to 1 and that $\delta_k(\omega^*,\PP_b)$ is positive for sufficiently large $k$.
\par For $x,y\ge1$, let $N(x,y;\N)\colonequals\#\{n \le x\colon\omega^*(n) \ge y\}$. Then we have
\[N(x,y;\N)=\sum_{k\ge y}\#(\mathcal{L}_k(\omega^*,\N)\cap[1,x]).\]
It has been shown \cite[Theorem 1]{FP} that there exists a suitable constant $c>0$ such that 
\[\left\lfloor\frac{x}{y^{c\log\log y}}\right\rfloor\le N(x,y;\N) \ll \frac{x\log y}{y}\]
for all $x\ge1$ and all sufficiently large $y$.
Here we consider the analogue of this counting function defined by $N(x,y;\PP_b)\colonequals\#\{b<p\le x\colon\omega^*(p-b) \ge y\}$ for any $b\in\Z\setminus\{0\}$. A simple adaptation of the proof of \cite[Theorem 1]{FP} yields the following theorem.
\begin{thm}\label{thm:N(x,y;P_b)}
For any $b\in\Z\setminus2\Z$, there exist constants $c_1,c_2>0$ depending on $b$, such that 
\[\frac{\pi(x)}{y^{c_1\log\log y}}< N(x,y;\PP_b) \ll \frac{\pi(x)\log y}{y}\]
for all sufficiently large $x$ and $y\le x^{c_2/\log\log x}$, where the implied constant in ``$\ll$" as well as the threshold for ``sufficiently large" may depend on $b$.
\end{thm}
Here the hypothesis $y\le x^{c_2/\log\log x}$ is required only for the lower bound to hold. The key ingredients in the proof of the upper bound are \cite[Theorem 3]{LMP}, which serves as a substitute for \cite[Theorem 1.2]{MPP} employed in the proof of \cite[Theorem 1]{FP}, and the estimate
\begin{align*}
\sum_{b<p\le x}\sum_{\substack{q-1\mid p-b\\ q-1\le z}}1\le\sum_{q\le z+1}\pi(x_b;q-1,b)&=\sum_{\substack{q\le z+1\\(q-1,b)=1}}\pi(x_b;q-1,b)+O(z)\\
&\ll\frac{x}{\log x}\sum_{\substack{q\le z+1\\(q-1,b)=1}}\frac{1}{\varphi(q-1)}+z\ll\pi(x)\log\log z
\end{align*}
for all $3\le z\le x$ with $\log z/\log x$ bounded away from 1, where we have used Brun--Titchmarsh and Proposition \ref{prop:shiftphi}. For the proof of the lower bound in Theorem \ref{thm:N(x,y;P_b)}, we would like to follow the proof of \cite[Theorem 1]{FP} as well. This is where the hypothesis $2\nmid b$ enters the picture. To achieve this, however, we need to upgrade \cite[Proposition 10]{APR} to the following result.
\begin{prop}\label{prop:maxomega*}
Let $b\in\Z\setminus2\Z$. There exist constants $c=c(b)>0$ and $x_0=x_0(b)\ge3$ with the property that for any integer $N>2$, one can find a square-free $M\in\N\cap[1,x^2)$ such that $\gcd(M,b)=1$, $N\nmid M$, and $\omega^*(M)\ge e^{c\log x/\log\log x}$ whenever $x\ge x_0$.
\end{prop}
Let us accept Proposition \ref{prop:maxomega*} for the moment and see how the lower bound in Theorem \ref{thm:N(x,y;P_b)} may be derived from it. We apply \cite[Proposition 8]{APR} with $\epsilon=1/2$, which gives us some $\delta\in(0,1)$ and $k_0(x)\in\N\cap((\log x)^{3/2},\infty)$ when $x$ is sufficiently large. By Proposition \ref{prop:maxomega*}, there exists an absolute constant $c>0$ such that we can find a square-free $M\in\N\cap[1,z^2)$ with $\gcd(M,b)=1$, $k_0(x)\nmid M$ and $\omega^*(M)\ge e^{c\log z/\log\log z}$ for all sufficiently large $z$. Take $z=y^{(c_1/3)\log\log y}$ with some suitable constant $c_1>0$ depending on $c$ and choose $c_2>0$ small enough depending on $c_1$ and $\delta$, so that we have $M<z^2\le x^{\delta}$ and $\omega^*(M)\ge y$ whenever $y\le x^{c_2/\log\log x}$ is sufficiently large. Moreover, according to \cite[Proposition 8]{APR}, we have
\[N(x,y;\PP_b)\ge\sum_{\substack{p\le x\\ p\equiv b\psmod{M}}}1\ge\frac{1}{\log x}\sum_{\substack{p\le x\\ p\equiv b\psmod{M}}}\log p>\frac{x}{2\varphi(M)\log x}>\frac{\pi(x)}{3z^2},\]
which yields the desired lower bound. Now we outline the proof of Proposition \ref{prop:maxomega*} below. 
\begin{proof}[Proof sketch of Proposition \ref{prop:maxomega*}]
Note that it suffices to consider the case where $b$ and $N$ are both square-free. We modify the proof of \cite[Proposition 10]{APR}. In what follows, we shall adopt the set-up and notation used there. Firstly, we fix a small $\epsilon>0$ whose value will be determined later (its original value is $1/4$ in the proof of \cite[Proposition 10]{APR}) and replace the original $k_1$ by $k_1'$ which is the product of primes $\le U\colonequals\max((1/4)\delta\log x,T)$ coprime to $b$, where $\delta\in(0,1/3)$ and $T>0$ originate from \cite[Proposition 9]{APR}. This will give us a new $k'$ in place of the original $k$. Following the proof of \cite[Proposition 10]{APR}, we have $k_0(x)\nmid (k'b)^2$. Next, we count the number $A(k',b,N)$ of $m\le x$ and $p\le x$ with $p-1$ square-free, satisfying
\[m(p-1)\equiv 0~\psmod{k'},~\gcd(m(p-1),b)=1\text{~and~}N\nmid m(p-1).\]
We define $A(k',b,\infty)$ to be the number of $m,p\le x$ satisfying the same conditions above except for $N\nmid m(p-1)$. The proof of \cite[Proposition 10]{APR} yields
\[A(k,1,\infty)>\frac{x^2}{20k\log x}\sum_{d\mid k}\frac{\varphi(d)}{d}\ge\frac{x^2}{k\log x}\left(\frac{3}{2}\right)^{\omega(k)}>\frac{x^2}{k\log x}\left(\frac{3}{2}\right)^{(1/4)\delta\log x/\log\log x}.\]
To complete the proof, it suffices to prove an inequality for $A(k',b,N)$ similar to this. It is not hard to see that the same argument yields
\begin{equation}\label{eq:A(k',b,infty)}
A(k',b,\infty)\gg\frac{x^2}{k'\log x}\sum_{d\mid k'}\frac{\varphi(d)}{d}\ge\frac{x^2}{k'\log x}\left(\frac{3}{2}\right)^{\omega(k')}>\frac{x^2}{k'\log x}\left(\frac{3}{2}\right)^{(1/6)\delta\log x/\log\log x},
\end{equation}
as long as $\epsilon$ is sufficiently small. This is where we need the assumption that $2\nmid b$, which ensures that 
\[\sum_{d\mid b}\mu(d)\frac{\psi(d)}{d}=\prod_{p\mid b}\frac{p(p-2)}{p^2-p-1}>0,\]
where the definition of $\psi(d)$ is given in \cite[Proposition 9]{APR}.
\par Suppose first that $N\ge x^{\delta/3}$. Then the number of $m\le x$ and $p\le x$ with $p-1$ square-free such that $m(p-1)\equiv 0~\psmod{N}$ is at most
\begin{align*}
\sum_{uv=N}\sum_{\substack{m\le x\\ u\mid m}}\sum_{\substack{p\le x\\ p\equiv 1\psmod{v}}}1\le\frac{2^{\omega(N)}}{N}x^2=o\left(\frac{x^2}{k'\log x}\left(\frac{3}{2}\right)^{(1/6)\delta\log x/\log\log x}\right) =o(A(k',b,\infty)).
\end{align*}
It follows from \eqref{eq:A(k',b,infty)} that
\begin{equation}\label{eq:A(k',b,N)}
A(k',b,N)\gg\frac{x^2}{k'\log x}\left(\frac{3}{2}\right)^{(1/6)\delta\log x/\log\log x}.
\end{equation}
\par Next, we suppose that $N<x^{\delta/3}$ has a prime factor $q\in(U,x^{\delta/3})$. Then we have $q\nmid k'$. For any $d\mid k'q$, the number of primes $p\le x$ with $p-1$ square-free and $d\mid p-1$ is 
\[\sum_{\substack{p\le x\\ p\equiv 1\psmod{d}}}\mu(p-1)^2=\sum_{n\le (x-1)/d}\mu(dn)^21_{\PP}(dn+1)=\sum_{n\le (x-1)/d}\mu(n)^21_{(n,d)=1}1_{\PP}(dn+1),\]
which, by Lemma \ref{lem:weightedHL} and the fact that $d\le\sqrt{x}$, is $\ll x/(d\log x)$. In addition, the number of $m\le x$ with $\gcd(m,k'q)=k'q/d$ is $\ll (x/k'q)\varphi(d)$. Thus, the number of $m\le x$ and $p\le x$ with $p-1$ square-free such that $m(p-1)\equiv0~\psmod{k'q}$ is 
\[\ll\frac{x^2}{k'q\log x}\sum_{d\mid k'q}\frac{\varphi(d)}{d}<\frac{2x^2}{k'q\log x}\sum_{d\mid k'}\frac{\varphi(d)}{d}=o\left(\frac{x^2}{k'\log x}\sum_{d\mid k'}\frac{\varphi(d)}{d}\right)=o(A(k',b,\infty)),\]
since $q>(1/4)\delta\log x$. Hence, we have \eqref{eq:A(k',b,N)} in this case as well.
\par Finally, we consider the case where $N<x^{\delta/3}$ and all of the prime factors of $N$ are $\le U$. If $\gcd(N,b)>1$, then $A(k',b,N)\ge A(k',b,\infty)$. So we may assume $\gcd(N,b)=1$. Then we have $N\mid k'$. Since $N>2$ is square-free, it has a prime factor $2<q\le U$. If we replace $k'$ by $k''=k'/q$, then the same argument in the proof of \cite[Proposition 10]{APR} still works. Hence, it suffices to prove a lower bound for $A(k'',b,N)$. It is clear that $A(k'',b,N)\ge A(k'',bq,\infty)$. A lower bound for $A(k'',bq,\infty)$ uniform in $2<q\le U$ that is analogous to \eqref{eq:A(k',b,infty)} can be established in the same fashion. This completes the proof.
\end{proof}

\vspace{2mm}
{\noindent\bf Acknowledgments.} The author thanks Carl Pomerance for helpful discussions and suggestions. He is also grateful to Paul Pollack for his constructive feedback on an earlier draft of the paper. Finally, he thanks the Max Planck Institute for Mathematics for its hospitality during his visit in 2024.

\end{document}